\documentclass[a4paper]{amsart}
\usepackage{amsthm,amsmath,amssymb,amscd,multirow}
\usepackage[latin1]{inputenc}
 \newtheorem{theorem}{Theorem}
 \newtheorem{proposition}[theorem]{Proposition}
 \newtheorem{lemma}[theorem]{Lemma}
 \newtheorem{corollary}[theorem]{Corollary}
 \newtheorem*{theoA}{Theorem A}
 \newtheorem*{theoB}{Theorem B}
 \newtheorem*{theoC}{Theorem C}
 \newtheorem*{theoD}{Theorem D}
 
 \theoremstyle{definition}
 \newtheorem{definition}[theorem]{Definition}
 \newtheorem{example}[theorem]{Example}

 \theoremstyle{remark}
 \newtheorem{remark}[theorem]{Remark}
 
 \numberwithin{equation}{section}
 \numberwithin{theorem}{section}

 \newcommand{\Spec}{{\rm Spec}}
 \newcommand{\Spf}{{\rm Spf}}
 \newcommand{\vol}{{\rm vol}}
 \newcommand{\Quot}{{\rm Quot}}
 \newcommand{\Stab}{{\rm Stab}}
 \newcommand{\End}{{\rm End}}
 \newcommand{\Nilp}{{\rm Nilp}}
 \newcommand{\Sets}{{\rm Sets}}
 \newcommand{\hei}{{\rm ht}}
 \newcommand{\id}{{\rm id}}
 
 \newcommand{\leng}{{\rm lg}}
 \newcommand{\et}{{\rm et}}
 \newcommand{\m}{{\rm m}}
 \newcommand{\bi}{{\rm bi}}
 
 \newcommand{\rk}{{\rm rk}}
 \newcommand{\determ}{{\rm det}}
 \newcommand{\Ann}{{\rm Ann}}
 \newcommand{\red}{{\rm red}}
 
 \newcommand{\defe}{{\rm def}}
 
  \newcommand{\Hom}{{\rm Hom}}
 \newcommand{\F}{{\boldsymbol F}}

\newenvironment{mylist}{\begin{list}{(\roman{enumi})}{\setlength{\leftmargin}{0pt}\setlength{\itemindent}{15pt} \usecounter{enumi}}}{\end{list}}

\begin{document}
\begin{title}
{Moduli spaces of $p$-divisible groups}
\end{title}
\author{Eva Viehmann}
\address{Mathematisches Institut der Universit\"{a}t Bonn\\ Beringstrasse 1\\53115 Bonn\\Germany}
\date{}
\begin{abstract}
We study the global structure of moduli spaces of
quasi-isogenies of $p$-divisible groups introduced by Rapoport and
Zink. We determine their dimensions and their sets of connected
components and of irreducible components. If the isocrystals of
the $p$-divisible groups are simple, we compute the cohomology of
the moduli space. As an application we determine which moduli
spaces are smooth.
\end{abstract}
\maketitle 

\section{Introduction}

Let $k$ be a perfect field of characteristic $p$ and $W=W(k)$ its
ring of Witt vectors. Let $\sigma$ be the Frobenius automorphism
on $k$ as well as on $W$. By $\Nilp_W$ we denote the category of
schemes $S$ over $\Spec(W)$ such that $p$ is locally nilpotent on
$S$. Let $\overline{S}$ be the closed subscheme of $S$ that is
defined by the ideal sheaf $p\mathcal{O}_S$. Let $\mathbb{X}$ be a
$p$-divisible group over $k$ with rational Dieudonn\'{e} module $(N,F)$. We assume that $\mathbb{X}$ is decent (see \cite{RapoportZink}, 2.14), i. e. $N$ is generated by elements $x$ satisfying an equation $F^sx=p^rx$ for some integers $s>0$ and $r$. Note that this condition is always satisfied if $k$ is algebraically closed.

We consider the functor
$$\mathcal{M}:\Nilp_W\rightarrow  \Sets,$$ which
assigns to $S\in \Nilp_W$ the set of isomorphism classes of pairs
$(X,\rho)$, where $X$ is a $p$-divisible group over $S$ and $\rho:
\mathbb{X}_{\overline{S}}=\mathbb{X}\times_{\Spec
(k)}\overline{S}\rightarrow X\times_{S}\overline{S}$ is a
quasi-isogeny. Two pairs $(X_1, \rho_1)$ and $(X_2, \rho_2)$ are
isomorphic if $\rho_1\circ\rho_2^{-1}$ lifts to an isomorphism
$X_2\rightarrow X_1$. This functor is representable by a formal
scheme $\mathcal{M}$, which is locally formally of finite type
over $\Spf(W)$ (see \cite{RapoportZink}, Thm. 2.16). Let
$\mathcal{M}_{\red}$ be its underlying reduced subscheme, that is the reduced subscheme of $\mathcal{M}$ defined by the maximal ideal of definition. The irreducible
components of $\mathcal{M}_{\red}$ are projective varieties over $\Spec (k)$ (\cite{RapoportZink}, Prop. 2.32).

These moduli spaces and their generalisations for moduli problems of type (EL) or (PEL) serve to analyse the local structure of Shimura varieties which have an interpretation as moduli spaces of abelian varieties. In \cite{RapoportZink} they are used to prove a uniformization theorem for Shimura varieties along Newton strata. Mantovan (see \cite{Mantovan}) computes the cohomology of certain (PEL) type Shimura varieties in terms of the cohomology of Igusa varieties and of the corresponding (PEL) type Rapoport-Zink spaces. In \cite{Fargues}, Fargues shows that the cohomology of basic unramified Rapoport-Zink spaces realises local Langlands correspondences.

For $p$-divisible groups whose rational Dieudonn\'{e} module is simple, the moduli spaces have been studied by de Jong and Oort in
\cite{deJongOort}. They show that the connected components are
irreducible and determine their dimension. In the general case
very little is known besides the existence theorem. This paper is
directed towards a better understanding of the global
structure of $\mathcal{M}_{\red}$.

We now state our main results. 

Let $\mathbb{X}=\mathbb{X}_{\m}\times \mathbb{X}_{\bi}\times
\mathbb{X}_{\et}$ be the decomposition of $\mathbb{X}$ into its
multiplicative, bi-infinitesimal, and \'{e}tale part. To formulate the
result about the set of connected components we exclude the
trivial case $\mathbb{X}_{\bi}=0$.

\begin{theoA} Let $\mathbb{X}$ be non-ordinary. Then $$\pi_0(\mathcal{M}_{\red})\cong \left(
GL_{\hei(\mathbb{X}_{\m})}(\mathbb{Q}_p)/GL_{\hei(\mathbb{X}_{\m})}(\mathbb{Z}_p)\right)\times \left(GL_{\hei(\mathbb{X}_{\et})}(\mathbb{Q}_p)/GL_{\hei(\mathbb{X}_{\et})}(\mathbb{Z}_p)\right)\times \mathbb{Z}$$ and the connected components are also rationally connected.
\end{theoA}

Next we consider the set of irreducible components of
$\mathcal{M}_{\red}$. From now on we assume that $k$ is algebraically closed. Let $(N,F)$ be the rational Dieudonn\'{e} module of $\mathbb{X}$ and $$J=\{g\in GL(N) \mid g\circ F=F\circ g\}.$$ There is an action of $J$ on $\mathcal{M}_{\red}$. We choose a
decomposition $N=\bigoplus_{j=1}^{l}N_{j}$ with $N_{j}$ simple of slope
$\lambda_j=m_j/(m_j+n_j)$ with $(m_j,n_j)=1$ and
$\lambda_j\leq\lambda_{j'}$ for $j<j'$. Let $M_0\subset N$ be the lattice associated to a minimal $p$-divisible group (see \cite{Oort} or Remark \ref{remm0}).

\begin{theoB}\begin{mylist}
\item The action of $J$ on the set of irreducible
components of $\mathcal{M}_{\red}$ is transitive and induces a bijection between this set and $J/(J\cap \Stab(M_0))$.
\item $\mathcal{M}_{\red}$ is equidimensional with
\begin{equation}\label{gldimform}
\dim\mathcal{M}_{\red}=\sum_{j}\frac{(m_j-1)(n_j-1)}{2}+\sum_{j<j'}m_jn_{j'}.
\end{equation}
\end{mylist}
\end{theoB}

In \cite{Oort3}, Oort defines an almost product structure (that is, up to a finite morphism) on Newton strata of moduli spaces of abelian varieties. It is given by an isogeny leaf and a central leaf for the $p$-divisible group. The dimension of the isogeny leaf is the same as that of the corresponding $\mathcal{M}_{\red}$. He announces a joint paper with Chai, in which they prove a dimension formula for central leaves (compare \cite{Oort3}, Remark 2.8). The dimension of the Newton polygon stratum itself is known from \cite{Oort1}. Then the dimension of $\mathcal{M}_{\red}$ can also be computed as the difference of the dimensions of the Newton polygon stratum and the central leaf.

Let $G= GL_h$ over $\mathbb{Q}_p$ where $h=\hei (\mathbb{X})$. Let $B$ be the Borel subgroup of lower triangular matrices. Let $\nu=(\lambda_1,\dotsc,\lambda_l)$ be the Newton vector associated to $(N,F)$. Here each $\lambda_j$ occurs $m_j+n_j$ times. Let $\mu=(0,\dotsc,0,1,\dotsc,1)$ with multiplicities $\dim\mathbb{X}$ and $h-\dim\mathbb{X}$. Let $\rho$ be the half-sum of the positive roots. Note that $J$ has an interpretation as the set of $\mathbb{Q}_p$-valued points of an algebraic group. Denote by $\rk_{\mathbb{Q}_p}$ the dimension of a maximal $\mathbb{Q}_p$-split subtorus. Let $\defe_G(F)=\rk_{\mathbb{Q}_p}G-\rk_{\mathbb{Q}_p}J$, which is equal to $h-l$. Then one can reformulate (\ref{gldimform}) as
\begin{equation}\label{gldimform2}
\dim \mathcal{M}_{\red}=\langle \rho,\mu-\nu\rangle -\frac{1}{2}\defe_G(F).
\end{equation}

Let $G$ be a split connected reductive group over a finite extension $\F$ of $\mathbb{Q}_p$ or $\mathbb{F}_p[[t]]$. Let $\mathcal{O}$ be the ring of integers in $\F$ and let $K=G(\mathcal{O})$. Let $L$ be the completion of the maximal unramified extension of $\F$ and let $\sigma$ be the Frobenius of $L$ over $\F$. Let $\mu$ be a conjugacy class of one-parameter subgroups of $G$ and $b\in B(G,\mu)$ (compare \cite{Rapoport}, 5). Let $\varepsilon^{\mu}$ be the image of $p$ or $t\in G(\F)$ under $\mu$. Let 
\begin{equation}
X_{\mu}(b)_K=\{g\in G(L)/K\mid g^{-1}b\sigma(g)\in K\varepsilon^{\mu} K\}
\end{equation}
be the generalised affine Deligne-Lusztig set associated to $\mu$ and $b$. There are two cases where it is known that $X_{\mu}(b)_K$ is the set of $\mathbb{F}$-valued points of a scheme, where $\mathbb{F}$ is the residue field of $\mathcal{O}_L$. The first case is that $\F=\mathbb{Q}_p$ and that $X_{\mu}(b)_k$ is the set of $k$-valued points of a Rapoport-Zink space of type (EL) or (PEL). In our case choose $\F=\mathbb{Q}_p$, $G=GL_h$ and choose an isomorphism $N\cong \Quot(W)^h$ such that $M_0$ corresponds to the lattice generated by the standard basis. Let $\mu$ be as above. We write $F=b\sigma$ with $b\in G$. Then a bijection is given by
\begin{eqnarray*}
X_{\mu}(b)_K&\rightarrow& \mathcal{M}_{\red}(k)\\
g&\mapsto&gM_0.
\end{eqnarray*}
The second case is that $\F$ is a function field. In this case, the dimension of the generalised affine Deligne-Lusztig variety has been determined in \cite{GHKR}, \cite{dimdlv}. The formula for the dimension is completely analogous to (\ref{gldimform2}). 

If the isocrystal of $\mathbb{X}_{\bi}$ is
not simple, Theorems A and B imply that the connected components of
the moduli space are not irreducible and thus not smooth. Now
assume that the isocrystal of $\mathbb{X}_{\bi}$ is simple of slope $m/(m+n)$. Then the
connected components are irreducible and projective. By
$\mathcal{M}_{\red}^0$ we denote the
connected component of the identity in the moduli space.

Let $m,n\in \mathbb{N}=\{0,1,\dotsc\}$ with $(m,n)=1$ be as above. A normalised
cycle is a $m+n$-tuple of integers $B=(b_0,\dotsc,b_{m+n-1})$ with
$b_0>b_i$, $b_{m+n-1}+m=b_0$, $\sum_i b_i=\sum_i i$ and
$b_{i+1}\in\{b_i+m,b_i-n\}$ for all $i$ (compare \cite{deJongOort}, 6). There are only finitely
many such cycles. Let $B^+=\{b_i\in B \mid b_{i}+m\in B\}$ and
$B^-=\{b_i\mid b_i-n\in B\}$. Then $B=B^+\sqcup B^-$. For $j\in
\mathbb{N}$ let $d(j)$ be the number of cycles $B$ such that
$\mathcal{V}(B)=\{(d,i) \mid b_d\in B^+, b_i\in B^-, b_i<b_d\}$
has $j$ elements.

\begin{theoC}
Let $\mathbb{X}$ be bi-infinitesimal and let its isocrystal be simple.
Let $m,n$ and $d(j)$ be as above. Let $l\neq p$ be prime. Then
$H^{2j+1}(\mathcal{M}_{\red}^0,\mathbb{Q}_l)=0$ for all $j$ and $H^{2j}(\mathcal{M}_{\red}^0,\mathbb{Q}_l)$ is a successive extension of $d(j)$ copies of $\mathbb{Q}_l(-j)$.
\end{theoC}

This description uses a paving of $\mathcal{M}^0_{\red}$ by affine spaces which resembles the description of the geometric points in
\cite{deJongOort}, 5. As an application we show

\begin{theoD} Let $\mathbb{X}$ be a $p$-divisible group over $k$. Then
$\mathcal{M}_{\red}^0$ is smooth if and only if one of the
following holds. Either $\dim\mathcal{M}_{\red}^0=0$ or the isocrystal
$N$ of $\mathbb{X}_{\bi}$ is simple of slope $2/5$ or
$3/5$. In this last case, $\mathcal{M}_{\red}^0\cong
\mathbb{P}^1$.
\end{theoD}
The condition $\dim\mathcal{M}_{\red}^0=0$ holds if and only if $\mathbb{X}$ is ordinary or the isocrystal of $\mathbb{X}_{\bi}$ is simple of slope $m/(m+n)$ with $\min\{m,n\}=1$.

The moduli spaces $\mathcal{M}$ are special cases of affine Deligne-Lusztig sets which have a scheme structure. It would be interesting to generalise our results to more general of these affine Deligne-Lusztig sets, in particular to other Rapoport-Zink spaces of type (EL) or (PEL).

\noindent{\it Acknowledgement.} I am grateful to M. Rapoport for
introducing me to these problems and for his interest and advice.
I thank Th. Zink for his interest in my work and T. Wedhorn for
many helpful discussions. Furthermore, I thank the referee for helpful comments.
\tableofcontents

\section{Review of methods}
Let $R$ be a commutative ring of characteristic $p>0$.
\subsection{Witt vectors}\label{secwitt}
Let $W(R)$ be the ring of Witt vectors of $R$. The Frobenius operator $\sigma:R\rightarrow R$ with $\alpha\mapsto \alpha^p$ induces an operator $W(R)\rightarrow W(R)$ which we also denote by $\sigma$. We will also write $a^{\sigma}$ instead of $\sigma(a)$.

Let $a\in R$. By $[a]=(a,0,\dotsc)\in W(R)$ we denote the
Teichm\"{u}ller representative of $a$. This defines a multiplicative embedding $R\rightarrow W(R)$.
\begin{remark}\label{remwitt}
Let $a=(a_0, a_1, \dotsc), b=(0,\dotsc,0,b_n,b_{n+1},\dotsc)\in W(R)$ and $\lambda\in R$. Then
\begin{eqnarray}
\label{glwittplus}a+b&=&(a_0,\dotsc,a_{n-1},a_n+b_n,c_{n+1},\dotsc)\\
\label{glwittmal}[\lambda]b&=&(0,\dotsc,0,\lambda^{p^n}b_n,d_{n+1},\dotsc)
\end{eqnarray}
with $c_i,d_i\in R$ for $i\geq n+1$. Assume that $b_n\in R^{\times}$ and $-a_nb_n^{-1}=\lambda^{p^n}$ for some $\lambda\in R$. Then from (\ref{glwittplus}) and
(\ref{glwittmal}) we get that
$$a+[\lambda] b=(a_0,\dotsc,a_{n-1},0,c_{n+1},\dotsc)$$ with
$c_i\in R.$
\end{remark}
Let $w_0:W(R)\rightarrow R$ be the Witt polynomial with $(a_0,a_1,\dotsc)\mapsto a_0$.

\subsection{Dieudonn\'{e} modules}
The Dieudonn\'{e} modules used in this paper to describe $p$-divisible groups over perfect fields are always the corresponding covariant Dieudonn\'{e} modules.

Let
\begin{equation}
\mathcal{D}(R)=W(R)[F,V]/(FV=VF=p,F\alpha=\alpha^{\sigma}F,\alpha V=V\alpha^{\sigma})
\end{equation} be the Dieudonn\'{e} ring of $R$.
Then each element $A\in\mathcal{D}(R)$ has a normal form as a finite sum
$$A=\sum_{i,j\geq 0}V^i x_{ij}F^j$$ with $x_{ij}\in W(R).$ If $R=K$ is a perfect field, $A$ can also be written as a possibly infinite sum
\begin{equation}\label{glforma}
A=\sum_{i,j\geq 0}[a_{ij}]V^iF^j
\end{equation}
with $a_{ij}\in K$ and $\mid i-j\mid$ bounded.

Let $M$ be a Dieudonn\'{e} module over a perfect field $k$ of characteristic $p$ and $N$ its rational Dieudonn\'{e} module. For a $k$-algebra $R$ denote
\begin{eqnarray*}
M_R&=&M\otimes_{W(k)}W(R)\\
N_R&=&N\otimes_{\Quot(W(k))}W(R)[1/p].
\end{eqnarray*}

A lattice in $N$ which is also a Dieudonn\'{e} module is called
a Dieudonn\'{e} lattice.

\subsection{Displays}\label{secdispl}

To fix notation we give a summary of some definitions and results of \cite{Zink2} on displays of $p$-divisible groups.

Let $R$ be an excellent $p$-adic ring and let $p$ be nilpotent in $R$.
\begin{definition}
A \emph{display over $R$} is a quadrupel $\mathcal{P}=(P,Q,F,V^{-1})$, where $P$ is a finitely generated projective $W(R)$-module, $Q\subseteq P$ is a submodule and $F$ and $V^{-1}$ are $\sigma$-linear maps, $F:P\rightarrow P$ and $V^{-1}:Q\rightarrow P$, such that the following properties are satisfied:
\begin{mylist}
\item Let $I_R$ be the ideal in $W(R)$ defined by the condition that the first Witt polynomial $w_0$ vanishes. Then $I_RP\subseteq Q\subseteq P$ and there exists a decomposition $P=L\oplus T$ into a direct sum of $W(R)$-modules such that $Q=L\oplus I_RT$. It is called a normal decomposition.
\item $V^{-1}:Q\rightarrow P$ is a $\sigma$-linear epimorphism.
\item For $x\in P$ and $w\in W(R)$ we have
$V^{-1}(\!\!~^Vwx)=wFx$ where $\!\!~^V\cdot:W(R)\rightarrow W(R)$
is the Verschiebung.
\end{mylist}
Besides, a nilpotence condition for $V$ is required, see
\cite{Zink2}, Def. 11.
\end{definition}

\begin{example}\label{ex13} If $M$ is the Dieudonn\'{e} module of a formal
$p$-divisible group $\mathbb{X}$ over a perfect field $k$, then
$(M,VM,F,V^{-1})$ is a display over $k$. We refer to it as the
display associated to $M$. In this case a normal decomposition is obtained as follows. We choose representatives $w_1,\dotsc,w_m$ in
$VM$ of a basis of the $k$-vector space $VM/pM$ and set $L=\langle
w_1,\dotsc, w_m\rangle_{W(k)}$. Similarly, we choose
representatives $v_1,\dotsc ,v_n$
of a basis of $M/VM$ and set $T=\langle v_1,\dotsc, v_n\rangle_{W(k)}$.
\end{example}

If $R$ is an excellent local ring or if $R/pR$ is of finite type over a field, there is an equivalence of categories between the category of displays over $R$ and the category of $p$-divisible formal groups over $\Spec (R)$ (\cite{Zink2}, Thm. 103).

To the base change of $p$-divisible groups corresponds a base change for displays. More precisely, let $S$ be another excellent ring and $\varphi:R\rightarrow S$ a morphism. Then for any display $\mathcal{P}=(P,Q,F,V^{-1})$ over $R$ there is an associated display $$\mathcal{P}_S=(P_S,Q_S,F_S,V_S^{-1})$$ over $S$ with $P_S=W(S)\otimes_{W(R)}P$, called the base change of $\mathcal{P}$ with respect to $\varphi$. We call the second component of the base change $Q_S$, although in general, we only have $Q_S\supseteq W(S)\otimes_{W(R)}Q$. For a definition of the base change see \cite{Zink2}, Def. 20.

\begin{definition}
An {\it isodisplay} over $R$ is a pair $(\mathcal{I},F)$ where
$\mathcal{I}$ is a finitely generated projective $W(R)\otimes
\mathbb{Q}\,$-module and $F:\mathcal{I}\rightarrow \mathcal{I}$ is a
$\sigma$-linear isomorphism.
\end{definition}
Let $\mathcal{P}=(P,Q,F,V^{-1})$ be a display over $R$. Then the pair $(P\otimes \mathbb{Q},F)$, where $F$ is the extension to $P\otimes\mathbb{Q}$, is an isodisplay over $R$.

Let $\mathbb{X}$ be a $p$-divisible group over $k$ and $N$ its rational Dieudonn\'{e} module. Let $R$ be a $k$-algebra of finite type and let
$\mathcal{P}=(P,Q,F,V^{-1})$ be a display over $R$ with $P\otimes
\mathbb{Q}\cong N_R$. Then by \cite{Zink2}, Prop. 66, this
isomorphism induces a quasi-isogeny between $\mathbb{X}_R$ and the
$p$-divisible group corresponding to $\mathcal{P}$.

\section{Connected Components}\label{secconncomp}

In this section we determine the set of connected components of
$\mathcal{M}_{\red}$. By $\mathcal{M}_{\red}^0$
we denote the connected component of
$(\mathbb{X}, \id)$ in $\mathcal{M}_{\red}$.

Let $\mathbb{X}=\mathbb{X}_{\m}\times \mathbb{X}_{\bi}\times
\mathbb{X}_{\et}$ be the decomposition of $\mathbb{X}$ into its
multiplicative, bi-infinitesimal, and \'{e}tale part. The moduli
spaces $\mathcal{M}(\mathbb{X}_{\m})$ and
$\mathcal{M}(\mathbb{X}_{\et})$ corresponding to $\mathbb{X}_{\m}$ and
$\mathbb{X}_{\et}$ are discrete. Assume that $\mathbb{X}$ is multiplicative or \'{e}tale. There is a basis of $N$ such that $F=p^{\alpha}\sigma$ with $\alpha\in \{0,1\}$. Let $K$ be a perfect field and a $k$-algebra. Recall that we set $h=\hei(\mathbb{X})$. Then
\begin{eqnarray*}
\mathcal{M}(\mathbb{X})(K)&\cong& \{g\in GL_h(\Quot(W(K)))/GL_h(W(K))\mid g^{-1}\sigma(g)\in GL_h(W(K))\}\\
&\cong& \{g\in GL_h(\Quot(W(K)))/GL_h(W(K))\mid g=\sigma(g)\}\\
&\cong& GL_h(\mathbb{Q}_p)/GL_h(\mathbb{Z}_p)
\end{eqnarray*}
independently of $K$. We define
\begin{equation}\Delta=\begin{cases}
\mathcal{M}(\mathbb{X}_{\m})\times \mathcal{M}(\mathbb{X}_{\et})\times
\mathbb{Z}& \text{if } \mathbb{X}_{\bi} \text{ is nontrivial}\\
\mathcal{M}(\mathbb{X}_{\m})\times\mathcal{M}(\mathbb{X}_{\et})&\text{else.}

\end{cases}\end{equation}
Let $S\in \Nilp_W$ and let $\rho:\mathbb{X}_{\overline{S}}\rightarrow X_{\overline{S}}$ be a
quasi-isogeny where $X$ is a $p$-divisible group over $S$. From
\cite{Messing}, Lemma II.4.8 we get a morphism
$X\twoheadrightarrow X_{\et}\rightarrow S$, such that $X_{\et}\rightarrow S$ is \'{e}tale, as well as a quasi-isogeny $\rho_{\et}:\mathbb{X}_{\et,\overline{S}}\rightarrow X_{\et,\overline{S}}$, functorially in $\rho$. This defines a morphism
$$\kappa_{\et}:\mathcal{M}_{\red}\rightarrow \mathcal{M}(\mathbb{X}_{\et}).$$ By duality one also obtains a
morphism $$\kappa_{\m}:\mathcal{M}_{\red}\rightarrow
\mathcal{M}(\mathbb{X}_{\m}).$$ Finally, the morphism
$\hei:\mathcal{M}_{\red}\rightarrow \mathbb{Z}$ maps a
quasi-isogeny to its height. Let
\begin{eqnarray}
\nonumber \kappa &:&
\mathcal{M}_{\red}\rightarrow \Delta\\
\kappa&=&\begin{cases}
(\kappa_{\m},\kappa_{\et},\hei)&\text{if } \mathbb{X}_{\bi} \text{ is nontrivial}\\
(\kappa_{\m},\kappa_{\et})&\text{else.}
\end{cases}
\end{eqnarray}
\begin{theorem}\label{thmconncomp}$\kappa$ identifies $\Delta$ with the set of rationally connected components of $\mathcal{M}_{\red}$. \end{theorem}
This follows from Lemma \ref{lemkappasurj} and Proposition \ref{lemfaserk}. The theorem implies that the action of $J$ on $\pi_0(\mathcal{M}_{\red})$ is transitive. Especially, all connected components are mutually isomorphic.
\begin{lemma}\label{lemkappasurj}
$\kappa$ is surjective.
\end{lemma}

\begin{proof}
It is enough to show that for nontrivial $\mathbb{X}_{\bi}$ there
is a quasi-isogeny $\mathbb{X}_{\bi,K}\rightarrow X$ of any given height over some algebraically closed
field $K$. There is a quasi-isogeny $\rho$ from
$\mathbb{X}_{\bi,K}$ to a product of groups that are up
to isogeny simple. Thus it is enough to show the statement for a bi-infinitesimal group that is simple up to isogeny. Especially, its dimension and its height are coprime. Let $M\subset N$ be the lattice corresponding to $(\mathbb{X},\id)$. Then the lattice $p^aV^b(M)$ corresponds to a $k$-valued point $(X,\rho)$ of $\mathcal{M}_{\red}$ with $\hei(\rho)=ah+b\cdot\dim \mathbb{X}$. Thus for every given integer we may choose $a$ and $b\in\mathbb{Z}$ such that $\hei(\rho)$ is equal to that integer.
\end{proof}

\begin{proposition}\label{lemfaserk}
Let $K$ over $k$ be a perfect field and $(X,\rho), (X',\rho')$ two
$K$-valued points of $\mathcal{M}_{\red}$ with
$\kappa((X,\rho))=\kappa((X',\rho'))$. Then the two points are in
the same connected component of $\mathcal{M}_{\red}$.
\end{proposition}
For the proof we need the following lemma.
\begin{lemma}\label{LemAustausch}
Assume that $\mathbb{X}$ is bi-infinitesimal. Let $M\subset N_K$ be a Dieudonn\'{e} lattice and
$$v_1, v_2 \in (F^{-1}M\cap V^{-1}M)\setminus M.$$ Then $\langle M,
v_1\rangle $ and $\langle M, v_2\rangle$ are Dieudonn\'{e}
lattices. There is a morphism $\mathbb{A}^1\rightarrow\mathcal{M}_{\red}$ mapping $0$ and $1$ to the $K$-valued points corresponding to $\langle M,
v_1\rangle $ and $\langle M, v_2\rangle$, respectively.
\end{lemma}
Here $\langle M, v\rangle$ denotes the $W(K)$-module generated by
$M$ and $v$.
\begin{proof}We may assume that $\langle M, v_1\rangle \neq \langle M,
v_2\rangle$. We define a quasi-isogeny of $p$-divisible groups
over $\Spec(K[t])$ such that $\langle M, v_1\rangle$ and $\langle
M, v_2\rangle$ are the lattices corresponding to the
specialisations at $t=1$ and $t=0$, respectively. To do this we
describe the corresponding subdisplay $(P',Q',F,V^{-1})$ of the
isodisplay $N_{K[t]}$ of $\mathbb{X}_{K[t]}$. We use the notation of
\ref{secdispl}. Let
\begin{eqnarray}
\label{gllocdeft}T&=&\langle Vv_1,Vv_2,w_1,\dotsc,w_{n-2}\rangle\\
\label{gllocdefl}L&=&\langle pv_1,pv_2,x_1,\dotsc,x_{m-2}\rangle
\end{eqnarray}
be a normal decomposition of the display associated to $M$. As the
classes of $Vv_1$ and $Vv_2$ in $M/VM$ are linearly independent
over $K$, we can choose such $w_i, x_i$ that the elements on
the right hand sides of (\ref{gllocdeft}) and (\ref{gllocdefl})
are representatives of bases of the $K$-vector spaces $T/pT$ and
$L/pL$. We now set
\begin{align*}
T'&=\langle [t]^{\sigma}\otimes v_1+[1-t]^{\sigma}\otimes v_2,1\otimes( Vv_1-Vv_2),1\otimes w_1,\dotsc,1\otimes w_{n-2}\rangle_{W(K[t])}\\
L'&=\langle [t]\otimes Vv_1+[1-t]\otimes Vv_2,
p\otimes(v_1-v_2),1\otimes x_1,\dotsc,1\otimes
x_{m-2}\rangle_{W(K[t])}
\end{align*}
and further $P'=L'+T'$ and $Q'=L'+I_{K[t]}T'$. Here $\langle \cdot
\rangle_{W(K[t])}$ denotes the $W(K[t])$-submodule of $N_{K[t]}$
generated by the elements in the brackets. We have 
$$\det \left(\begin{array}{cccc}
p^{-1}[t]^{\sigma}&1&0&0\\
p^{-1}[1-t]^{\sigma}&-1&0&0\\
0&0&[t]&1\\
0&0&[1-t]&-1
\end{array}\right)=p^{-1}([t]+[1-t])^{\sigma}([t]+[1-t]).$$ As $W(K[t])$ is complete with respect to its augmentation ideal, $[t]+[1-t]$ is a unit. Hence, the determinant is a unit in $W(K[t])\otimes \mathbb{Q}$. Thus $P'\otimes \mathbb{Q}=N_{K[t]}$ and $L'\cap T'=(0)$. To show that $P'$ and $Q'$ define a subdisplay it remains to verify that $V^{-1}$ is a $\sigma$-linear epimorphism
from $Q'$ to $P'$. This follows from $$V^{-1}([t]\otimes
Vv_1+[1-t]\otimes Vv_2)=[t]^{\sigma}\otimes v_1 +
[1-t]^{\sigma}\otimes v_2.$$

The specialisations of this display for
$t=0$ and $t=1$ are as desired.
\end{proof}

\begin{proof}[Proof of Proposition \ref{lemfaserk}]

As $K$ is a perfect field, we can decompose $X$ and $\rho$ into
$$\rho=(\rho_{\m},\rho_{\bi},\rho_{\et}):\mathbb{X}_{\m,K}\times \mathbb{X}_{\bi,K}
\times \mathbb{X}_{\et,K}\rightarrow X_{\m}\times X_{\bi}\times
X_{\et},$$ and similarly for $\rho'$. The morphism $\kappa$ maps
$\rho$ to $(\rho_{\m},\rho_{\et}, \hei(\rho))$. The assumption implies
$\rho_{\m}=\rho'_{\m}$, $\rho_{\et}=\rho'_{\et}$ and
$\hei(\rho)=\hei(\rho')$. Assume that the proposition is proved
for $\mathbb{X}_{\bi}$. Then we can construct a quasi-isogeny over
a connected base $S=\overline{S}$ with fibres $\rho$ and $\rho'$ by extending
a quasi-isogeny with fibres $\rho_{\bi}$ and $\rho'_{\bi}$ on
$\mathbb{X}_{\bi,S}$ by the constant isogeny
$(\rho_{\m}\times\rho_{\et})_S=(\rho'_{\m}\times\rho'_{\et})_S$ on
$(\mathbb{X}_{\m}\times\mathbb{X}_{\et})_S$. Thus for the rest of the
proof we may assume that $\mathbb{X}=\mathbb{X}_{\bi}$.

From the two quasi-isogenies we get Dieudonn\'{e} lattices $M,
M'\subset N_K$ with $\vol(M)=\vol(M')$. We prove the proposition
by induction on the length of $M'/M\cap M'$. If the length is 0,
the lattices are equal and the statement is trivial. Let now
$M'\neq M$. As $\mathbb{X}$ is bi-infinitesimal, both $F$ and $V$ are topologically nilpotent on $M$. As $M\cap M'\subsetneq M$, there is an element
\begin{equation}v_1\in
M\setminus (FM+VM+M').
\end{equation}
Let further $v'\in M'\setminus M$. Let $i'$ be maximal with
$F^{i'}v'\notin M$ and $j'$ maximal with $V^{j'}F^{i'}v'=v_2\notin
M$. Then
\begin{equation} v_2\in M'\cap F^{-1}(M'\cap M)\cap V^{-1}(M'\cap
M)\setminus M.
\end{equation}
Let $\{v_1,x_1,\dotsc x_l\}$ be a basis of the $K$-vector space
$M/(FM+VM+(M'\cap M)).$ We choose representatives of the $x_i$ in
$M$, which we also denote by $x_i$. Let $\widetilde{M}$ be the
smallest $\mathcal{D}(K)$-module containing $FM, VM, M'\cap M$,
and all $x_i$. By the choice of the $x_i$ we have $v_1\notin
\widetilde{M}$. As $\widetilde{M}\subset M$ we have $v_2\notin
\widetilde{M}$. We also get $Fv_2,Vv_2 \in M'\cap M\subseteq
\widetilde{M}$. Thus the tuple $(\widetilde{M},v_1,v_2)$ satisfies
the assumption of Lemma \ref{LemAustausch}. Hence there is a morphism from the affine line to the moduli space whose image contains $\langle
\widetilde{M},v_1 \rangle=M$ and $\langle
\widetilde{M},v_2\rangle$. As $M'\cap M\subseteq \widetilde{M}$ and $v_2\in M'\setminus \widetilde{M}$, the length of $M'/(\langle
\widetilde{M},v_2\rangle\cap M')$ is smaller than that of
$M'/(M'\cap M)$. Thus the assertion
follows from the induction hypothesis.
\end{proof}


\section{Irreducible Components}\label{secirrcomp}
From now on we assume that $k$ is algebraically closed.
\subsection{Statement of the Theorem}\label{secirr1}
To formulate the main theorem of this section we need some
notation. We introduce a system of generators for the rational Dieudonn\'{e} module $N$ of $\mathbb{X}$. Let
\begin{equation}
N=\bigoplus_{j=1}^{j_0}N_{\lambda_j}
\end{equation}
be the isotypic decomposition of $N$ with $N_{\lambda_j}$ isoclinic of slope $\lambda_j$ and $\lambda_j<\lambda_{j'}$
for $j<j'$. There are coprime integers $0\leq m_j\leq h_j$ with $h_j>0$ and
$\lambda_j=m_j/h_j$. Let $n_j=h_j-m_j$. Note that this convention is the opposite of the one made in \cite{deJongOort}, meaning that we use $n$ instead of $m$ and vice versa. For
each $j$ we choose $a_j, b_j\in \mathbb{Z}$ with 
\begin{equation}\label{glab}a_jh_j+b_jm_j=1.\end{equation} We
define additive maps $\pi_j:N\rightarrow N$ by
\begin{equation}\label{gldefpi}
\pi_{j}\mid_{N_{\lambda_{j'}}}=\begin{cases}p^{a_j}F^{b_j}& \text{if }j'=j\\
\id_{N_{\lambda_{j'}}}&\text{else}
\end{cases}\end{equation}
and
$\sigma_{j}:N\rightarrow N$ by
\begin{equation}
\sigma_{j}\mid_{N_{\lambda_{j'}}}=\begin{cases}V^{-m_j}F^{n_j}& \text{if }j'=j\\
\id_{N_{\lambda_{j'}}}&\text{else.}
\end{cases}\end{equation}

There is an algebraic group $J=J_N$ over $\mathbb{Q}_p$ associated to the moduli problem and the isocrystal $N$, see \cite{RapoportZink}, 1.12. For each $\mathbb{Q}_p$-algebra $R$ its $R$-valued points are defined as $$J_N(R)=\{g\in GL(N\otimes_{\mathbb{Q}_p}R) \mid g\circ F=F\circ g\}.$$ In the following we will write $J$ or $J_N$ instead of $J_N(\mathbb{Q}_p)$ to simplify the notation.

\begin{remark}\label{remj}
Let $g\in GL(N)$. Then $g\in J_{N}$ if and only if $g$ commutes
with all $\pi_j$ and $\sigma_j$. Indeed, $g\in J_{N}$ if and only
if $g=\bigoplus_{j}g|_{N_{\lambda_j}}$ and
$g|_{N_{\lambda_j}}\in J_{N_{\lambda_j}}$ for all $j$. On
$N_{\lambda_j}$ we have $\pi_j=p^{a_j}F^{b_j}$ and
$\sigma_j=p^{-m_j}F^{m_j+n_j}$, and for the other direction
$F=\pi_j^{m_j}\sigma_j^{a_j}$.
\end{remark}

Let
\begin{equation}
N_{\lambda_j}=\bigoplus_{i=1}^{l_j}N_{j,i}
\end{equation}
be a decomposition into simple isocrystals. Let $e_{ji0}\in
N_{j,i}\setminus\{0\}$ with
\begin{equation}\label{gldefe0}F^{h_j}e_{ji0}=p^{m_j}e_{ji0}.\end{equation}
For $l\in\mathbb{Z}$ let
\begin{equation}e_{jil}=\pi_j^{l}e_{ji0}.\end{equation}
By (\ref{gldefe0}) the $e_{jil}$ are independent of the choice of
$a_j$ and $b_j$ in (\ref{glab}). Besides,
\begin{equation}\begin{aligned}
e_{j,i,l+h_j}&=\pi_{j}^{l+h_j}(e_{ji0})\\
&=\pi_j^lp^{h_ja_j}F^{h_jb_j}e_{ji0}\\
&=\pi_j^lp^{1-m_jb_j}F^{h_jb_j}e_{ji0}\\
\label{glpundpi}&=pe_{jil}
\end{aligned}\end{equation}
and analogously
\begin{align}F(e_{jil})&=e_{j,i,l+m_j},\\
V(e_{jil})&=e_{j,i,l+n_j},\\
\sigma_{j'}(e_{jil})&=e_{jil}
\end{align} for $1\leq j'\leq j_0$.
The $e_{jil}$ with $0\leq l<h_j$ form a basis of $N_{j,i}$ over
$\Quot (W(k))$.

Let $K\supseteq k$ be a perfect field. For $a\in K$ let $[a]\in W(K)$
be the Teichm\"{u}ller representative as in \ref{secwitt}. By
(\ref{glpundpi}) each $v\in N_K$ can be written as
\begin{equation}
v=\sum_{j=1}^{j_0}\sum_{i=1}^{l_j}\sum_{l\in\mathbb{Z}}[a_{jil}]e_{jil}
\end{equation}
with $a_{jil}\in K$ and $a_{jil}=0$ for $l$ small enough.

\begin{definition}\label{defm0}\begin{mylist}
\item Let $M_0\subset N$ be the lattice generated by the $e_{jil}$ with $l\geq 0$. 

\item For a lattice $M$ in some sub-isocrystal
$\widetilde{N}\subseteq N$ let
\begin{equation}
\vol_{\widetilde{N}}(M)=\leng((M_0\cap\widetilde{N})/(M_0\cap
M))-\leng(M/(M_0\cap M)).
\end{equation}
If $\widetilde{N}=N$ we write $\vol$ instead of $\vol_N$.
\end{mylist}
\end{definition}
\begin{remark}\label{remm0}
\begin{mylist}
\item $\pi_j(M_0)\subseteq M_0$ and $\sigma_j(M_0)\subseteq M_0$ for
all $j$, and $\vol(M_0)=0$.
\item We have $M_0=\bigoplus_{i,j} M_0\cap N_{j,i}$ and $\End(M_0\cap N_{j,i})$ is a maximal order in $\End(N_{j,i})$. Thus $M_0$ is the lattice associated to a minimal $p$-divisible group in the sense of Oort, compare \cite{Oort}. Besides, $J\cap \Stab (M_0)$ is the set of units in a maximal order of the semi-simple algebra $J=\End(N,F)$. 
\end{mylist}
\end{remark}

Using this notation we can formulate the main result of this section.
\begin{theorem}\label{thmirrcomp}
\begin{mylist}
\item There is a bijection between the set of irreducible
components of $\mathcal{M}_{\red}$ and $J/(J\cap \Stab(M_0))$.
\item $\mathcal{M}_{\red}$ is equidimensional with
\begin{equation*}
\dim\mathcal{M}_{\red}=\sum_{j}l_j\frac{(m_j-1)(n_j-1)}{2}+\sum_{(j,i)<(j',i')}m_jn_{j'},
\end{equation*}
where the pairs $(j,i)$ are ordered lexicographically.
\end{mylist}
\end{theorem}
The theorem shows that the action of $J$ on the set of irreducible components is transitive. Hence all irreducible components are isomorphic.

To prove the theorem, we first define an open and dense subscheme $\mathcal{S}_1\subseteq\mathcal{M}_{\red}$. It is enough to prove the corresponding statements for $\mathcal{S}_1$. We attach to each lattice $M$ corresponding to an element of $\mathcal{S}_1(k)$ the smallest lattice containing $M$ which corresponds to a minimal $p$-divisible group. These lattices form a single $J$-orbit, and the map leads to a morphism $\mathcal{S}_1\rightarrow J/(J\cap \Stab (M_0))$. We show that the inverse image in $\mathcal{S}_1$ of each element of $J$ is irreducible and of the dimension claimed in the theorem.

We assume until Section \ref{secmultet} that $\mathbb{X}$ is bi-infinitesimal. The general case is discussed in Section \ref{secmultet}.

\subsection{Definition of the subscheme $\mathcal{S}_1$}

\begin{definition}
Let $\mathcal{S}_1\subseteq\mathcal{M}_{\red}$ be the open
subscheme defined by the following condition. An $R$-valued point $(X,\rho)$ of $\mathcal{M}_{\red}$ lies in $\mathcal{S}_1$ if the display $(P,Q,F,V^{-1})$ over $R$ of $X$ has the property that the $R$-module $P/(Q+F(P))$ is locally free of rank 1. Here, $(Q+F(P))$ denotes the $W(R)$-submodule of $P$ generated by $F(P)$ and $Q$.
\end{definition}

\begin{remark}\label{lem35}
As we assumed $\mathbb{X}$ to be bi-infinitesimal, the rank of $P/(Q+F(P))$ is always positive. Thus $\mathcal{S}_1$ is open. 

Let $K \supseteq k$ be a perfect field and let $M\subset N_K$ be the lattice
associated to $x\in\mathcal{M}_{\red}(K)$. Recall that the $a$-invariant $a(M)$ of a Dieudonn\'{e} lattice $M$ over $K$ is defined as the dimension of the $K$-vector space $M/(FM+VM)$. Thus the point $x$ lies in $\mathcal{S}_1$ if and only if $a(M)=1$. Assume that this is the case. As $F$ and $V$ are topologically nilpotent on $M$, there is a $v\in M$ such that $M$ is the $\mathcal{D}(K)$-submodule of $N_K$ generated by $v$. Note that the $a$-invariant can also be defined without using the Dieudonn\'{e} module: If $X$ is a $p$-divisible group with Dieudonn\'{e} module $M$, then $a(M)=a(X)$ where $a(X)=\dim\Hom_K(\alpha_p,X)$. 
\end{remark}
\begin{lemma}
The open subscheme $\mathcal{S}_1\subseteq \mathcal{M}_{\red}$ is dense.
\end{lemma}
\begin{proof} Let $X_0$ be the $p$-divisible group of a $K$-valued point in $\mathcal{M}_{\red}\setminus \mathcal{S}_1$. By Proposition 2.8 of \cite{Oort1}, there exists a deformation of $X_0$ with constant Newton polygon such that the $a$-invariant at the generic fibre is 1. By \cite{OortZink}, Cor. 3.2 we get a corresponding deformation of the quasi-isogeny after a suitable base change preserving the generic fibre.
\end{proof}

\subsection{$K$-valued points of $\mathcal{S}_1$}\label{seckpoints}
Let $K \supseteq k$ be a perfect field. In this
section we classify the $K$-valued points of $\mathcal{S}_1$ by
introducing a normal form for the corresponding lattices in $N_K$.
We will write $N_{j,i}$ instead of $(N_{j,i})_K$.

\begin{lemma}\label{lemzerlv}
Let $M\subset N_K$ be the lattice associated to a $K$-valued point
of $\mathcal{S}_1$ and $v$ a generator of $M$ as $\mathcal{D}(K)$-submodule of $N_K$ as in Remark \ref{lem35}. Let $g$ be such that $v_p(\determ (g))$ is maximal among the $g\in J$ with $M\subseteq gM_0$. Then
\begin{equation}\label{gldarstv}
v=\sum_{j=1}^{j_0}\sum_{i=1}^{l_j}\sum_{l\geq
0}[a_{jil}]g(e_{jil})
\end{equation} with $a_{jil}\in K$ and
\begin{equation}\label{condstar}
\text{for each $j$, the $a_{ji0}$ for $1\leq i\leq l_j$ are linearly independent over }
\mathbb{F}_{p^{h_j}}.
\end{equation}
\end{lemma}

\begin{proof}We may assume that $N$ is isoclinic, as otherwise we can write
$v$ as a sum of elements of the $(N_{\lambda_j})_K$ and show the claim
for each summand separately. Assume that there is a nontrivial relation
\begin{equation}\label{glrelv}
\sum_{i=1}^{l_1}\alpha_i a_{1i0}=0
\end{equation}
 with $\alpha_i\in
\mathbb{F}_{p^{h_1}}$. After permuting the simple summands of $N_K$ we may assume that $\alpha_1$ is nonzero. Then we may also assume that $\alpha_1=-1$. We define $\delta\in GL(N_K)$ by
\begin{equation}\label{gldefdelta}
\delta(e_{1il})=\begin{cases}e_{1,1,l+1}&\text{if }i=1\\
e_{1il}+[\alpha_i]^{\sigma^{lb_1}} e_{1,1,l}&\text{if }i\geq
2\end{cases}
\end{equation}
for $l\in\mathbb{Z}$. This map is well defined as
\begin{eqnarray*}\delta(pe_{1il})&=&\delta(e_{1,i,l+h_1})\\
&=&pe_{1il}+[\alpha_i]^{\sigma^{(l+h_1)b_1}}e_{1,1,l+h_1}\\
&=&pe_{1il}+[\alpha_i]^{\sigma^{lb_1}} e_{1,1,l+h_1}\\
&=&p\delta(e_{1il})\end{eqnarray*} for $i>1$. We
also have
\begin{eqnarray*}
\delta\circ F(e_{1il})&=&\delta(e_{1,i,l+m_1})\\
&=&e_{1,i,l+m_1}+[\alpha_i]^{\sigma^{(l+m_1)b_1}}e_{1,1,l+m_1}\\
&=&e_{1,i,l+m_1}+[\alpha_i]^{\sigma^{lb_1+1-a_1h_1}}e_{1,1,l+m_1}\\
&=&F(e_{1il}+[\alpha_i]^{\sigma^{lb_1}}e_{1,1,l})\\
&=&F\circ\delta(e_{1il}),
\end{eqnarray*} for $i>1$ and $\delta\circ F(e_{11l})=F\circ\delta(e_{11l})=e_{1,1,l+m_1+1}$, so $\delta\in J$. Besides, $v_p(\determ(\delta))=1$. (\ref{gldarstv}) and (\ref{glrelv}) imply that $v\in g\circ\delta(M_0)$. As $v$ generates $M$, we have $M\subseteq g\circ \delta(M_0)$ in contradiction to the maximality of $v_p(\determ(g))$.
\end{proof}

\begin{remark}
Let $g$ be as in the lemma and $h\in \Stab(M_0)\cap J$. Then $g\circ h$ again satisfies the conditions of Lemma \ref{lemzerlv}.
\end{remark}

Let $M\subset N_K$ be a Dieudonn\'{e} lattice. Let $P(M)$ be the smallest
$\mathcal{D}(K)$-submodule of $N_K$ containing $M$ with
$$\sigma_j(P(M))\subseteq P(M)$$ and $$\pi_j (P(M))\subseteq P(M)$$
for all $j$. There exists a $c\in \mathbb{Z}$ with $M\subseteq
p^{c}M_0$. As $\pi_j(M_0)\subseteq M_0$ and $\sigma_j(M_0)\subseteq
M_0$ for all $j$, we get $M\subseteq P(M)\subseteq p^cM_0$. Hence
$P(M)$ is also a lattice in $N_K$. As all $\pi_j$ and $\sigma_j$
commute with $J$, we have $P(gM)=gP(M)$ for all $g\in
J$.

\begin{lemma}\label{lemtypm} Let $M$ be the Dieudonn\'{e} lattice corresponding to a
$K$-valued point of $\mathcal{S}_1$ and let $g$ be as in Lemma \ref{lemzerlv}. Then
\begin{equation}P(M)=gM_0.\end{equation} Especially, the class of $g$ in
$J/(\Stab(M_0)\cap J)$ is uniquely determined by $M$.
\end{lemma}

\begin{proof}
The inclusion $P(M)\subseteq gM_0$ follows from $v\in gM_0$ and $P(gM_0)=gM_0$. As all $\pi_j$ and
$\sigma_j$ commute with $J$ we may assume $g=\id$. On $N_{\lambda_j}$, the map $\pi_j$ is elementwise topologically nilpotent, while it is the
identity on each $N_{\lambda_{j'}}$ with $j'\neq j$.  Thus $v\in P(M)$
implies that its image $v_j=\sum_{i,l}[a_{jil}]e_{jil}$ under the projection to $(N_{\lambda_j})_K$ is in $P(M)$ for each $j$. Thus we may assume that $N=N_{\lambda_1}$ is isoclinic. Then $\pi_1$ is topologically nilpotent on $M_0$. Thus from Nakayama's lemma we see that $M_0\subseteq P(M)+\pi_1(M_0)$ implies the desired inclusion $M_0\subseteq P(M)$. To show $M_0\subseteq P(M)+\pi_1(M_0)$ we consider the image of $P(M)$ in $M_0/\pi_1(M_0)\cong K^{l_1}$.

Note that the set of solutions of a nontrivial equation of the form $\sum_{l=0}^{l_1-1}c_lx^{\sigma^{h_1l}}$ with $c_l\in K$ is an $\mathbb{F}_{p^{h_1}}$-vector space. The degree of the polynomial shows that its dimension is less or equal $l_1-1$. As the $l_1$ elements $\sigma_1^{h_1l}((a_{1i0})_i)$ for $0\leq l \leq l_1-1$ are linearly independent over $\mathbb{F}^{p^{h_1}}$, they cannot satisfy a common equation of this form. Hence they are linearly independent over $K$ as elements of $M_0/\pi_1M_0$. In other words, $v,\sigma_1(v),\dotsc,\sigma_1^{l_1-1}(v)$ is a basis of $M_0/\pi_1(M_0)$, and $M_0\subseteq P(M)+\pi_1(M_0)$.
\end{proof}

\begin{theorem}\label{lemvolm}
Let
\begin{equation} \label{glvorv}
v=\sum_{j=1}^{j_0}\sum_{i=1}^{l_j}\sum_{l\geq
0}[a_{jil}]g(e_{jil})\in N_K
\end{equation} with $a_{jil}\in K$ satisfying (\ref{condstar}) and let $g\in J$. Let $M$ the smallest $\mathcal{D}(K)$-submodule of $N_K$
containing $v$. Then
\begin{mylist}
\item $M$ is a lattice in $N_K$.
\item The volume of $M$ is of the form \begin{equation}
\vol(M)=v_p(\determ(g))+c,
\end{equation}
with
\begin{equation}c=\sum_{j}l_j\frac{(m_j-1)(n_j-1)}{2}+\sum_{(j,i)<(j',i')}m_jn_{j'},
\end{equation}
where the pairs $(j,i)$ are ordered lexicographically.
\item Let $I=I(N)\subset \coprod_{j,i}\mathbb{N}$ with
\begin{equation}\label{gldefI}
(\coprod_{j,i}\mathbb{N})\setminus
I=\{(j,i,l) \mid \,l=am_j+bn_j+\sum_{(j',i')<(j,i)}m_{j'}n_{j} \text{ for
some }a,b\geq 0\}.
\end{equation}
Then $ \mid I \mid =c$ and $(1,1,0)\notin I$, but $(j,i,0)\in I$ for $(j,i)\neq (1,1)$. There is an element $w\in M$ such that
\begin{equation}w=\sum_{(j,i)}\sum_{l\geq 0}[b_{jil}]g(e_{jil})
\end{equation}
with (\ref{condstar}) for $b_{jil}$ instead of $a_{jil}$. Besides, $b_{110}=1$, and $b_{jil}=0$ if
$(1,1,0)\neq(j,i,l)\notin I$. This element $w$ is a generator of
$M$ as $\mathcal{D}(K)$-module, and is called a normalised
generator. It only depends on the choice of the representative
$g\in [g]\in J/(J\cap \Stab(M_0))$.
\item Let
\begin{equation}\label{gllocelemm}x=\sum_{j,i,l}[b_{jil}]g(e_{jil})\in M.
\end{equation} We order the set $\coprod_{j,i}\mathbb{N}$ lexicographically. Then the index of the first nonzero summand of $x$ is not in $I$.
\item For all $(j,i,l)\notin I$ there is an $x\in M$ as in (\ref{gllocelemm}) such that $[b_{jil}]$ is its first nonzero coefficient.
\end{mylist}
\end{theorem}
For the proof of the theorem we need two technical
lemmas.
\begin{lemma}\label{lemabsch}
Let $a,b,m,n\in \mathbb{N}$ with $an+bm>mn$ and
$\lambda=m'/(m'+n')\in (0,1)$. Then
$$an'+bm'>\min\{nm',mn'\}.$$
\end{lemma}
\begin{proof}We may assume $0<a\leq m$ and $b>0$, because otherwise
the implication is evident. If the claim were false, then $\lambda(n+a-b)\geq a$ and
$\lambda(m+b-a)\leq m-a$. Especially, this implies $n+a-b>0$. From
our assumptions we get $m+b-a>0$. Thus $$a(m+b-a)\leq\lambda
(m+b-a)(n+a-b)\leq (m-a)(n+a-b),$$ in contradiction to $an+bm>mn$.
\end{proof}

\begin{lemma}\label{lemconsta}Let $v$ and $M$ be as in the theorem and assume $g=\id$ and $a_{110}=1$. Then there is an $A\in \mathcal{D}(K)$ of the form
\begin{equation*}A=F^{n_1}-V^{m_1}+\sum_{k > m_1n_1}[\alpha_{k}]V^{a(k)}F^{b(k)}\end{equation*}
with $\alpha_k\in K$ and the following properties:
\begin{mylist}
\item For $k>m_1n_1$, the exponents $a(k),b(k)$ are the unique positive integers with $-n_1< a(k)-b(k)\leq m_1$ and $a(k)n_1+b(k)m_1=k$.
\item $Av\in N'=\bigoplus_{(j,i)\neq(1,1)} N_{j,i}$.
\item $Av$ generates $M\cap N'$ as a Dieudonn\'{e}
lattice in $N'$.
\item Let $g'\in J_{N'}$ with $$g'(e_{jil})=e_{j,i,l+m_1n_j}.$$ Then
\begin{equation*} Av=\sum_{(j,i)\neq (1,1)}\sum_{l\geq
0}[a'_{jil}]g'(e_{jil}),
\end{equation*} with $a'_{jil}\in K$ such that for each $j$ all $a'_{ji0}$ (with
$1\leq i\leq l_j$ if $j\neq 1$ and $2\leq i\leq l_1$ if $j=1$) are
linearly independent over $\mathbb{F}_{p^{h_j}}$.
\end{mylist}
\end{lemma}

\begin{proof}
As $a_{110}=1$ and $F^{n_1}(e_{110})=V^{m_1}(e_{110})=e_{1,1,m_1n_1}$, we have
\begin{equation}
F^{n_1}v-V^{m_1}v=\sum_{l>m_1n_1}[c_{11l}]e_{11l}+\sum_{(j,i)\neq(1,1)}\sum_{l\geq
0}[c_{jil}]e_{jil}
\end{equation}
with $c_{jil}\in K$ for all $j,i,l$. For each $k>m_1n_1$ choose
$a(k)$ and $b(k)$ as in (i). Then
\begin{eqnarray*}
V^{a(k)}F^{b(k)}v&=&
\sum_{j,i}\sum_{l\geq 0}[a_{jil}^{\sigma^{b(k)-a(k)}}]e_{j,i,l+a(k)n_j+b(k)m_j}\\
&=&e_{11k}+\sum_{l>k}[d_{11l}]e_{11l}+\sum_{(j,i)\neq(1,1)}\sum_{l\geq0}[d_{jil}]e_{jil}
\end{eqnarray*} with $d_{jil}\in K$.
By Remark \ref{remwitt} we can inductively define $\alpha_k\in K$ for $k>m_1n_1$ such that 
\begin{equation*}
(F^{n_1}-V^{m_1}+\sum_{k=m_1n_1+1}^{k_0}[\alpha_k]V^{a(k)}F^{b(k)})(v)
\end{equation*}
is a linear combination of the $e_{11l}$ with $l>k_0$ and of an element of $N'$. Thus the element
\begin{equation}\label{glwahla}
A=F^{n_1}-V^{m_1}+\sum_{k>m_1n_1}[\alpha_k]V^{a(k)}F^{b(k)}
\end{equation}
satisfies $Av\in N'$. As $a(k)-b(k)$ is bounded, we have $A\in\mathcal{D}(K)$.

Let $B=\sum_{a,b\geq 0}[\beta_{ab}]V^aF^b\in \mathcal{D}(K)$ with
$\beta_{ab}\in K$ and $Bv\in M\cap N'$. We want to show that $B=CA$ for some $C\in \mathcal{D}(K)$. We assume $B\neq 0$. For
each $k\in \mathbb{N}$ such that there exists a $\beta_{ab}\neq 0$
with $an_1+bm_1=k$ let
$$d_k(B)=\min\{a-b\mid an_1+bm_1=k , \beta_{ab}\neq 0\}$$
and
$$d^k(B)=\max\{a-b\mid an_1+bm_1=k , \beta_{ab}\neq 0\}.$$ Furthermore let
\begin{eqnarray}
d(B)&=&\min\{d_k(B)\}\\
\tilde{d}(B)&=&\max\{d^k(B)\}
\end{eqnarray} The
existence of this minimum and maximum is equivalent to
$B\in\mathcal{D}(K)$. Inductively we construct $C_{k}\in
\mathcal{D}(K)$ with the following properties:
\begin{mylist}
\item[(a)] The coefficient of $V^cF^d$ in the representation of $B-\sum_{k'\leq k}C_{k'}A$ as in (\ref{glforma}) vanishes for all $c,d$ with $cn_1+dm_1\leq k$.
\item[(b)] If there exists a $\beta_{ab}\neq 0$ with $an_1+bm_1=k$, then $d(C_k)\geq d(B)$ and
$\tilde{d}(C_k)\leq \tilde{d}(B)$. Otherwise, $C_k=0.$
\item[(c)] If $B-\sum_{k'\leq k}C_{k'}A\neq 0$ then
\begin{eqnarray*}
d(B)&\leq&d(B-\sum_{k'\leq k}C_{k'}A)\\
\tilde{d}(B)&\geq&\tilde{d}(B-\sum_{k'\leq k}C_{k'}A).
\end{eqnarray*}
\end{mylist}
If $C=\sum_{k\geq 0}C_{k}$ exists in $\mathcal{D}(K)$, then this implies $B=CA$. By
replacing $B$ by $B-\sum_{k'< k}C_{k'}A$ we may assume that $k$ is
the least integer such that there exist $a,b$ with
$an_1+bm_1=k$ and $\beta_{ab}\neq 0$. We want to show that $d_k(B)\neq d^k(B)$. Assume that $d_k(B)=d^k(B)$. Then there is only one $\beta_{a_0b_0}\neq 0$ with $a_0n_1+b_0m_1=k$. Denote by $p_1$ the projection to $N_{1,1}$. We have
\begin{align}
\nonumber 0&=B(p_1(v))\\
\nonumber
&=[\beta_{a_0b_0}]V^{a_0}F^{b_0}(p_1(v))+\sum_{\{(a,b)\mid an_1+bm_1>k\}}[\beta_{ab}]V^{a}F^{b}(p_1(v))\\
\nonumber &=\sum_{l\geq
0}[\beta_{a_0b_0}a_{11l}^{\sigma^{b_0-a_0}}]e_{1,1,l+k}+
\sum_{\{(a,b)\mid an_1+bm_1>k\}}\sum_{l\geq
0}[\beta_{ab}a_{11l}^{\sigma^{b-a}}]e_{1,1,l+an_1+bm_1}.
\end{align}
Hence the coefficient of $e_{1,1,k}$ in the expression above is
$[\beta_{a_0b_0}a_{110}^{\sigma^{b_0-a_0}}]=[\beta_{a_0b_0}]$. This implies
$\beta_{a_0b_0}=0,$ a contradiction. Thus $d^k(B)>d_k(B)$. Note
that $m_1+n_1$ divides $d^k(B)-d_k(B)$. Let
$a,b$ with $a-b=d^k(B)$ be the pair of indices realising the maximum. Let $C_{k,1}=[-\beta_{ab}]V^{a-m_1}F^{b}$. From $d^k(B)\geq d_k(B)+m_1+n_1$ we see that $d_k(B)<a-b-m_1<d^k(B)$ and that $d^k(B)>d^k(B-C_{k,1}A)$ and
$d_k(B)\leq d_k(B-C_{k,1}A)$. Hence $d^k(B)-d_k(B)>d^k(B-C_{k,1}A)-d_k(B-C_{k,1}A)$. Using
a second induction on this difference, we can construct $C_k$ as a
finite sum of such expressions $C_{k,1}$. The fact that each pair $(a,b)$ occurs at most once in the construction of some $C_k$ together with (b) implies that the sum $C=\sum_{k\geq 0}C_k$ exists in $\mathcal{D}(K)$. This proves (iii).
 
Now we want to show (iv). We have
\begin{eqnarray*}
Av&=& F^{n_1}v-V^{m_1}v+\sum_{k>m_1n_1}[\alpha_k]V^{a(k)}F^{b(k)}v\\
&=&\sum_{(j,i)\neq(1,1)}\sum_{l\geq 0}\left(
F^{n_1}[a_{jil}]e_{jil}-V^{m_1}[a_{jil}]e_{jil}+
\sum_{k>m_1n_1}[\alpha_k]V^{a(k)}F^{b(k)}[a_{jil}]e_{jil}\right)\\
&=&\sum_{(j,i)\neq(1,1)}\sum_{l\geq 0}\left(
[a_{jil}^{\sigma^{n_1}}]e_{j,i,l+n_1m_j}-[a_{jil}^{\sigma^{-m_1}}]e_{j,i,l+m_1n_j}\right.\\
&&\hspace{5cm}+\sum_{k>m_1n_1}\left.[\alpha_k
a_{jil}^{\sigma^{b(k)-a(k)}}]e_{j,i,l+a(k)n_j+b(k)m_j}\right).
\end{eqnarray*}
For each $j$ and $i$ we determine the first nonvanishing
coefficient of some $e_{jil}$. First we consider summands with
$j=1$ and $i>1$. In this case $V^{a(k)}F^{b(k)}e_{1i0}=e_{1ik}$
with $k>m_1n_1$. Thus a candidate for the first coefficient is
that of $e_{1,i,m_1n_1}$, namely
$[a_{1i0}^{\sigma^{n_1}}-a_{1i0}^{\sigma^{-m_1}}]$. (Here we used
Remark \ref{remwitt} to determine
$w_0([a_{1i0}^{\sigma^{n_1}}]-[a_{1i0}^{\sigma^{-m_1}}])$.) As in
the proof of Lemma \ref{lemtypm} one sees that these coefficients
are again linearly independent over $\mathbb{F}_{p^{h_1}}$. Now we
consider summands with $j>1$. From Lemma \ref{lemabsch} and the
ordering of the $\lambda_j$ we get
$$a(k)n_j+b(k)m_j>\min\{n_1m_j,m_1n_j\}=m_1n_j.$$ Thus the first
nonzero coefficient is that of $e_{j,i,m_1n_j}$, namely
$a_{ji0}^{\sigma^{-m_1}}$. For fixed $j$ the $a_{ji0}$ were
linearly independent over $\mathbb{F}_{p^{h_j}}$, hence the new
first nonzero coefficients are again linearly independent. This
proves (iv).
\end{proof}
\begin{corollary}
Let $N$ be bi-infinitesimal and simple and $v\in
N_K\setminus\{0\}$. Then $\Ann(v)\subset \mathcal{D}(K)$ is a
principal left ideal.
\end{corollary}

\begin{proof}[Proof of Theorem \ref{lemvolm}]
Both $F$ and $V$ commute with $g$. Thus $M=gM'$ where $M'$ is generated by $$g^{-1}v=\sum_{j,i}\sum_{l\geq 0}[a_{jil}]e_{jil}.$$ Hence we may assume that $g=\id$.

Assertion (iii) is implied by (iv) and (v). We show that (i) and (ii) also follow from (iv) and (v). We consider the $\mathcal{D}(K)$-modules
\begin{equation}M^{j_0i_0l_0}=\langle M, \{e_{jil} \mid (j,i,l)\geq (j_0,i_0,l_0), l\geq 0\}\rangle_{W(K)}.
\end{equation}
Then $$M^{110}=M_0.$$ Using (v) one sees
$$M^{j_0,l_{j_0},d}=M$$ where
$d=\sum_{(j,i)\neq(j_0,l_{j_0})}m_jn_{j_0}+(m_{j_0}-1)(n_{j_0}-1)$.
For $(j,i,l)<(j',i',l')$ we have $M^{jil}\supseteq M^{j'i'l'}$
with equality if and only if
$$I\cap \{(j_1,i_1,l_1) \mid (j,i,l)\leq(j_1,i_1,l_1)<(j',i',l')\}=\emptyset.$$ Indeed, by (iv) and (v) this is equivalent to the condition
that  for each $(j_1,i_1,l_1)$ with
$(j,i,l)\leq(j_1,i_1,l_1)<(j',i',l')$, there is already an element
of $M$ whose first nonzero coefficient has index $(j_1,i_1,l_1)$.
As the colength of $M^{j,i,l+1}$ in $M^{jil}$ is at most $1$, this
implies $\vol(M)=\mid I\mid =c$.

We now prove (iv) and (v) using induction on $\sum_{j=1}^{j_0}l_j$, the number of simple summands of $N_K$.

By multiplying $v$ with $[a_{110}^{-1}]\in W(K)^{\times}$ we may assume that
$a_{110}=1$. First we consider the case that $N$ is simple. We have
$$V^aF^bv=e_{1,1,an_1+bm_1}+\sum_{l>0}[a_{11l}^{\sigma^{b-a}}]e_{1,1,l+an_1+bm_1}.$$
Thus for all $l$ that can be written as $l=an_1+bm_1$ with
$a,b\geq 0$ there is an element in $M$ of the form
$e_{11l}+\sum_{l'>l}[b_{l'}]e_{11l'}$, proving (v). Let now $x\in
M\setminus\{0\}$ and assume that $$x=\sum_{l\geq
0}[b_{11l}]e_{11l}$$ such that (iv) is not satisfied for $x$. Let
$[b_{11l_0}]$ with $(1,1,l_0)\in I$ be the first nonvanishing coefficient. We also have
a representation $x=\sum_{a,b\geq 0}[c_{a,b}]V^aF^b (v)$. Let
$(a_0,b_0)$ be a pair with $c_{a_0,b_0}\neq 0$ and minimal
$a_0n_1+b_0m_1$. Then $a_0n_1+b_0m_1\leq l_0$. As no
$l>(m_1-1)(n_1-1)$ is in $I$, we get $a_0n_1+b_0m_1\leq
(m_1-1)(n_1-1)$. Especially, $(a_0,b_0)$ is the unique pair of
nonnegative integers $(a,b)$ with $an_1+bm_1=a_0n_1+b_0m_1$. Hence
the coefficient of $e_{1,1,a_0n_1+b_0m_1}$ is the first nonzero
coefficient of $x$, proving (iv).

Let now $N$ be the sum of more than one simple summand. Let
$p_1:N_K\rightarrow N_{1,1}$ be the projection and
$$N'=\bigoplus_{(j,i)\neq(1,1)} N_{j,i}.$$
Note that $p_1(M)$ is the lattice in $N_{1,1}$ generated by
$p_1(v)$. Thus the theorem applied to the simple isocrystal $N_{1,1}$ yields that each $x\in
M\setminus (M\cap N')$ satisfies (iv), and that for each $(1,1,l)\notin
I$ there is an element $x\in M$ as in (v). We now consider elements
of $M\cap N'$. Let $I(N')$ be the index set corresponding to $N'$ as in (\ref{gldefI}), viewed as a subset of $\coprod_{(j,i)\neq (1,1)}\mathbb{N}$. Then
one easily checks that
\begin{equation}\label{glrecI}
I\cap\coprod_{(j,i)\neq (1,1)}\mathbb{N}= \{(j,i,l) \mid (j,i,l-m_1n_j)\in I(N')\}.
\end{equation}
In Lemma \ref{lemconsta} we proved that $M\cap N'$ is generated by
$Av$ for some $A\in \mathcal{D}(K)$ and determined the
corresponding $g'\in J_{N'}$, which only shifts the last
indices of the basis by $m_1n_j$. Thus the induction hypothesis
implies that there is an $x\in M\cap N'$ with $x=\sum_{(j,i)\neq
(1,1)} [c_{jil}]e_{jil}$ and first nonzero coefficient $[c_{jil}]$
if and only if $(j,i,l-m_1n_j)\notin I(N')$. Together with
(\ref{glrecI}), this implies the theorem.
\end{proof}

\subsection{Irreducible subvarieties of $\mathcal{S}_1$}\label{secirrvar}

Let $I$ be the index set defined in Theorem \ref{lemvolm}(iii).
Denote the coordinates of a point in
$\mathbb{A}_{k}^{I}$ by $a_{jil}$ with $(j,i,l)\in I$. 
Let $U=U(N)\subseteq \mathbb{A}_k^I$ be the affine open
subvariety defined by (\ref{condstar}). Let $a_{110}=1$. Then $U$ is defined by the condition that for each $j$, the $a_{ji0}$ for
$1\leq i\leq l_j$ have to be linearly independent over
$\mathbb{F}_{p^{h_j}}$. We write $U=\Spec(R)$.

For each $g\in J$ we want to define a morphism
$$\varphi_g:U\rightarrow \mathcal{S}_1.$$ For $g=\id$ we describe the
corresponding quasi-isogeny of $p$-divisible groups over $U$ via
the display of the $p$-divisible group. As $J$ acts on
$\mathcal{M}_{\red}$ we can define $\varphi_g$ for general $g$ as
the composition of $\varphi_{\id}$ and the action of $g$.

We may choose $\mathbb{X}$ inside its isogeny class in such a way that its Dieudonn\'{e} module is $M_0$. Let $(P,Q,F,V^{-1})$ be the base change of the display of $\mathbb{X}$ from $k$ to $R$. Let
$$v=e_{110}+\sum_{(j,i,l)\in I}[\sigma^{\sum_{j,i}m_j}(a_{jil})]e_{jil}\in P$$ and
\begin{eqnarray}
\label{glglobT}\tilde{T}&=&\langle v,Fv, \dotsc, F^{\sum_{i,j}n_j-1}v\rangle_{W(R)}\\
\label{glglobL}\tilde{L}&=&\langle Vv,\dotsc, V^{\sum_{i,j}m_j}v\rangle_{W(R)}
\end{eqnarray}
as $W(R)$-submodules of $P$. Let
$\tilde{P}=\tilde{L}+\tilde{T}$ and
$\tilde{Q}=I_R\tilde{T}+\tilde{L}$. We have to show that
$(\tilde{P},\tilde{Q},F,V^{-1})$ is a display where $F$ and
$V^{-1}$ are the restrictions of $F$ and $V^{-1}$ on $P$. By construction we have $I_R\tilde{P}\subset
\tilde{Q}\subset \tilde{P}$ and $\tilde{P}$ and $\tilde{Q}$ are
finitely generated $W(R)$-modules. The results of the preceding section show that the reduction of $P/\tilde{P}$ in a $K$-valued point of $U$ is torsion of some finite index which is bounded by a constant only depending on $N$ and not on the specific point. Thus there is a power of $p$ annulating $P/\tilde{P}$, and $\tilde{P}\otimes
\mathbb{Q}\cong P\otimes\mathbb{Q}$. For dimension reasons,
$\tilde{P}$ has to be free and $\tilde{L}\cap
\tilde{T}=(0)$. Hence $\tilde{L}$ and $\tilde{T}$ form a normal
decomposition. The third condition for a display and
the nilpotence condition on $V$ are satisfied because they were
satisfied on $P$. We now determine the matrix associated to $F|_{\tilde{T}}$
and $V^{-1}|_{\tilde{L}}$ as in \cite{Zink2},(9) to show that the image of $F$
is again in $\tilde{P}$ and that $V^{-1}:\tilde{Q}\rightarrow \tilde{P}$
is a $\sigma$-linear epimorphism. The matrix  is of the following form:
$$\left(\begin{array}{ccccccccc}
0&\cdots&0&\ast&~&1\\
1&\ddots&\vdots&\vdots\\
~&\ddots&0&\vdots&~&&&\multirow{2}{6mm}[10pt]{\LARGE 0}\\
&~&1&\ast\\
~\\
&~&~&\ast&~&0&1\\
~&\multirow{2}{3mm}{\LARGE 0}&&\vdots&&~&\ddots&\ddots\\
~&&~&\vdots&~&\multirow{2}{0mm}[5pt]{\LARGE 0}&~&\ddots&1\\
~&&~&\ast&&~&~&~&0
\end{array}
\right)$$
All columns except the one corresponding to $F^{\sum_{j,i}n_j-1}v$ have exactly one nonzero entry, which is 1. We now have to show that the remaining column has entries in $W(R)$. We use induction on $\sum l_j$, the number of simple summands of $N$, to show the following property. Let $v=\sum_{j,i,l\geq 0}[b_{jil}]e_{jil}\in N_R$ for some $R$ such that for each $j$ all non-trivial linear combinations of its coefficients $b_{ji0}$ with coefficients in $\mathbb{F}_{p^{h_j}}$ are in $R^{\times}$ and that $b_{110}=1$. Assume furthermore that the coefficients of $v$ are in $\sigma^{\sum_{j,i}m_j}(R)$. Then $$F^{\sum_{j,i}n_j}v=\sum_{0<k<\sum_{j,i}n_j} \gamma_{-k}F^{k}v+\sum_{0\leq k\leq \sum_{j,i}m_j}\gamma_kV^kv$$ with $\gamma_k\in W(R)$.

Let $A\in\Ann(p_1(v))\subseteq\mathcal{D}(R)$ of the same form as in Lemma \ref{lemconsta}. The construction for this over $R$ is the same as over $K$. As we chose $a(k)-b(k)\leq m_1$, all coefficients of $(F^{n_1}-V^{m_1})(v)$ and $V^{a(k)}F^{b(k)}(v)$ are in $\sigma^{\sum_{(j,i)\neq(1,1)}m_j}(R)$. Thus the coefficients of $A$ are also in $\sigma^{\sum_{(j,i)\neq(1,1)}m_j}(R)$. We can write $A$ in the form
\begin{equation}\label{gldarsta}
A=\sum_{0<k\leq n_1}\alpha_{-k} F^{k}+\sum_{0\leq k\leq m_1}\alpha_k V^k
\end{equation}
with $\alpha_k\in W(\sigma^{\sum_{(j,i)\neq(1,1)}m_j}(R))$ and $\alpha_{-n_1}=1$. The last equation holds as $-n_1<a(k)-b(k)$ for all $k$. Hence if $N$ is simple, the equation $Av=0$ gives the desired relation for $F^{n_1}v$. We now consider the case that $N$ is not simple. As in the proof of Lemma \ref{lemtypm}, the linear independence condition on the coefficients of $v$ implies the following similar condition for $Av$. For each $j$, all non-trivial linear combinations with coefficients in $\mathbb{F}_{p^{h_j}}$ of the first coefficients of the projections of $Av$ on all $(N_{j,i})_R\subseteq N'_R$ are invertible in $R$. Especially the projection of $Av$ on the second simple summand of $N_R$ is nonzero, and its first nonzero coefficient $[\beta]$ is invertible. This implies that $[\beta^{-1}]Av$ is (up to an index shift as in Lemma \ref{lemconsta}(iv)) an element of $N'_R$ satisfying the conditions needed to apply the induction hypothesis. Thus \begin{multline*}F^{\sum_{(j,i)\neq(1,1)}n_j}([\beta^{-1}]Av)\\=\sum_{0<k<\sum_{(j,i)\neq(1,1)}n_j} \gamma_{-k}F^{k}([\beta^{-1}]Av)+\sum_{0\leq k\leq \sum_{(j,i)\neq(1,1)}m_j}\gamma_kV^k([\beta^{-1}]Av)
\end{multline*} with $\gamma_k\in R$. Together with (\ref{gldarsta}) this leads to the desired relation for $F^{\sum_{j,i}n_j}v$. Hence $F^{\sum_{j,i}n_j}v\in \tilde{P}$ and $(\tilde{P},\tilde{Q},F,V^{-1})$ is a display. 

This display, together with the identity as
isomorphism of isodisplays, induces a quasi-isogeny of
$p$-divisible groups over $U$, that is a morphism
$\varphi_{\id}:U\rightarrow \mathcal{M}_{\red}$. For all $g\in J$ let
$$\mathcal{S}(g)=\varphi_{g}(U)=g\circ \varphi_{\id}(U).$$ As can be seen on $K$-valued points, the subvarieties $\mathcal{S}(g)$ and $\mathcal{S}(g')$ for $g,g'\in J$ are equal if and only if $[g]=[g']$ in $J/(J\cap \Stab(M_0))$.
\begin{remark}
Let $M\subset N$ be a lattice and let $\mathcal{P}=(P,Q,F,V^{-1})$ be the display associated to an $S$-valued point of $\mathcal{M}_{\red}$. We consider $P$ as a submodule of $N_S$. Then the condition that $P$ is contained
in $M_S$ is a closed condition on $S$.
\end{remark}
\begin{lemma}\label{lemende}
For all $g\in J/(J\cap \Stab(M_0))$, the subscheme
$\mathcal{S}(g)$ is a connected component of $\mathcal{S}_{1}$.
\end{lemma}
\begin{proof} Let $B$ be a set of representatives of $J/(J\cap \Stab(M_0))$. Then the $\mathcal{S}(g)$ for $g\in B$ are disjoint and cover
$\mathcal{S}_1$. This holds as it is true for their sets of
$K$-valued points for every algebraically closed $K$, compare Lemma \ref{lemtypm}. The height
of the quasi-isogeny is constant on each connected component of
$\mathcal{M}_{\red}$ and thus of $\mathcal{S}_1$. Let $M$ be a
lattice associated to a $K$-valued point of $\mathcal{S}_1$. Then
$\vol(M)-\vol(P(M))=c$ is a constant only depending on $N$. If $M\subseteq gM_0$ for some $g\in J$ with $v_p(\determ (g))=\vol(P(M))$, then $P(M)=gM_0$. Thus $\mathcal{S}(g)(K)$ consists of the lattices $M$ with $\vol(M)=v_p(\determ(g))+c$ and $M\subseteq gM_0$. Hence $\mathcal{S}(g)$ is closed. The fact that $\mathcal{M}$ is locally
formally of finite type implies that the disjoint union is locally
finite. Thus the $\mathcal{S}(g)$ are also open.
\end{proof}

\subsection{The general case}\label{secmultet}
Now we consider the case of general $\mathbb{X}$ over $k$, that is
we do not assume that $\mathbb{X}$ is bi-infinitesimal. The results obtained for the set of irreducible
components and the dimension in the bi-infinites\-imal case also hold
in this more general context. To see this, we again consider the
set of $K$-valued points for an algebraically closed field $K$.
Over $K$ each quasi-isogeny $\rho:\mathbb{X}_K\rightarrow X$
splits into a product of quasi-isogenies between the \'{e}tale,
multiplicative, and bi-infinitesimal parts of $\mathbb{X}_K$ and
$X$. The results of Section \ref{secconncomp} show that the
connected component of the point $x\in\mathcal{M}_{\red}(K)$
corresponding to $\rho$ is given by fixing the \'{e}tale and
multiplicative part of the quasi-isogeny and its height. Thus all
points of one connected component may be classified by considering
the bi-infinitesimal parts of the quasi-isogenies. Quasi-isogenies
corresponding to the irreducible subvarieties of
$\mathcal{M}_{\red}$ of Section \ref{secirrvar} can be defined in
this context as a product of a constant quasi-isogeny of the
\'{e}tale and multiplicative parts of $\mathbb{X}_K$ and the
quasi-isogeny of Section \ref{secirrvar} for the bi-infinitesimal part.


\section{Cohomology}
\subsection{A paving of $\mathcal{M}_{\red}$ for $N$ simple}
Let $\mathbb{X}$ be a bi-infinitesimal $p$-divisible group over an
algebraically closed field $k$ of characteristic $p$ whose
rational Dieudonn\'{e} module $N$ is simple. In this case we obtain from \cite{deJongOort}, Thm. 5.11 that $\mathcal{M}^0_{\red}$ is irreducible. Hence
it is also projective (compare \cite{RapoportZink},
Prop. 2.32).

In the following we simplify the notation by writing $e_l$ instead of $e_{1,1,l}$ for the basis of $N$ and $n=n_1$, $m=m_1$, and $h=m+n$. We pave $\mathcal{M}_{\red}^0$ with affine spaces to compute its cohomology. This is inspired by a description of
the geometric points of $\mathcal{M}_{\red}$ in \cite{deJongOort}. Note however that we interchanged the roles of $n$ and $m$ with respect to the convention in de Jong and Oort's paper.

Let $K$ be a perfect field over $k$. We recall a combinatorial
invariant for $K$-valued points of $\mathcal{M}_{\red}$ from
\cite{deJongOort}, 5. A subset $A\subseteq \mathbb{Z}$ is called a {\it semimodule} if it is bounded below and satisfies $m+A\subseteq A$ and $n+A\subseteq A$. It is called {\it normalised} if $ \mid \mathbb{N}\setminus A(M) \mid
= \mid A(M)\setminus \mathbb{N} \mid $. One easily sees that there are only
finitely many normalised semimodules. In fact, their number is
$\binom{h}{m}/h$, see \cite{deJongOort}, 6.3. For every semimodule $A$, there is a unique integer $l$ such that $l+A$ is normalised. We call $l+A$ the normalisation of $A$. Each element of $N_K$ can be
uniquely written as $\sum_l [a_l]e_{l}$ with $a_l\in K$ and
$a_l=0$ for $l$ small enough. We call the least $l\in
\mathbb{Z}$ with $a_l\neq 0$ the {\it first index} of the element. Let
$M\subset N_K$ be the lattice associated to $x\in \mathcal{M}_{\red}(K)$.
As $M$ is a Dieudonn\'{e} lattice,
\begin{equation}A=A(M)=\{l\in \mathbb{Z} \mid l \text{ first index of some }v\in M\}
\end{equation}
is a semimodule called the semimodule of $x$ or $M$. From the definition of the volume we get
\begin{equation}
\vol(M)= \mid \mathbb{N}\setminus A(M) \mid - \mid A(M)\setminus \mathbb{N} \mid .
\end{equation}
We may assume that $\id_{\mathbb{X}}$ corresponds to a lattice of
volume $0$. Then the semimodules of $K$-valued points of
$\mathcal{M}_{\red}^0$ are normalised. 

\begin{proposition}
For each normalised semimodule $A$ there is a locally closed reduced
subscheme $\mathcal{M}_A\subseteq \mathcal{M}_{\red}^0$ which is defined
by the property that for each perfect field $K$, the set
$\mathcal{M}_A(K)$ consists of the points with semimodule $A$. The
$\mathcal{M}_A$ are disjoint and cover $\mathcal{M}_{\red}^0$.
\end{proposition}
\begin{proof}
It is enough to show that for every normalised semimodule $A$
there is an open subscheme $\mathcal{M}_{\leq A}$ of
$\mathcal{M}_{\red}^0$, such that for every perfect field $K$ the
set $\mathcal{M}_{\leq A}(K)$ consists of all points whose
semimodules $A'$ satisfy the following condition: There is a
bijection $f:A'\rightarrow A$ with $f(a)\geq a$ for all $a$. A
Dieudonn\'{e} lattice $M$ of volume $0$ corresponds to an element
of $\mathcal{M}_{\leq A}(K)$ if and only if for all $a\in A$, the
length of $M/\langle e_{a+1},e_{a+2},\dotsc\rangle_{W(K)}$ is at
least $ \mid A\cap\mathbb{Z}_{\leq a} \mid $. Every normalised semimodule $A'$ contains an element $a_0\leq 0$. All $a> mn-m-n\geq a_0+mn-m-n$ can be written as $a=a_0+\alpha m+\beta n$ with $\alpha,\beta\geq 0$. Thus $a\in A'$ for all $a> mn-m-n$ and all normalised semimodules $A'$. Hence the condition above is an intersection of finitely many open conditions on $\mathcal{M}_{\red}^0$.
\end{proof}

We want to identify each $\mathcal{M}_A$ with an affine space. To
do this we need further combinatorial invariants from
\cite{deJongOort}. Let $A$ be a semimodule. We arrange the $h$ elements of $A\setminus
(h+A)$ in the following way. Let $b_0$ be the largest element.
For $i=1,\dotsc,h-1$ we choose inductively $b_i\in A\setminus
(h+A)$ to be $b_{i-1}-n$ or $b_{i-1}+m$, depending on which of
the elements lies in $A\setminus (h+A)$. Then $b_0=b_{h-1}+m$.
The tuple
\begin{equation}B=B(A)=(b_0,\dotsc,b_{h-1})
\end{equation} is
called the cycle of $A$. One can recover $A$ as $A=\{b_i+lh
\mid b_i\in B, l\geq 0\}$. This defines a bijection between the
set of semimodules and the set of cycles, that is of $h$-tuples
of integers $b_i$ satisfying $b_0>b_i$, $b_{h-1}+m=b_0$ and
$b_i\in \{b_{i-1}-n,b_{i-1}+m\}$ for all $i\neq 0$. The
normalisation condition for semimodules is equivalent to the
condition
\begin{equation}\sum_i
b_i=\frac{(h-1)h}{2}.\end{equation} We split each cycle $B$
in two parts:
\begin{eqnarray}
B^+&=&\{b_i \mid b_i+m\in B\}\\
B^-&=&\{b_i \mid b_i-n\in B\}.
\end{eqnarray}
Let \begin{equation}\mathcal{V}(B)=\{(d,i) \mid b_d\in B^+, b_i\in B^-,
b_i<b_d\}\end{equation} and \begin{equation}R=k[a_{d,i} \mid (d,i)\in
\mathcal{V}(B)],\end{equation}
\begin{equation}S=\Spec(R)=\mathbb{A}_k^{\mathcal{V}(B)}.\end{equation} We
define a quasi-isogeny $\mathbb{X}_S\rightarrow X$ by describing
the display of $X$ as a subdisplay of the isodisplay
$N_R$ of $\mathbb{X}_S$. For each $b_i\in B$ we
want to define an element $v_i\in N_R$ which has first
index $b_i$ and first coefficient 1 in all closed points of $S$. We want
the $v_i$ to satisfy the following relations:
\begin{equation}\label{glstart} v_0=e_{b_0}
\end{equation} and
\begin{equation}\label{gldefdisp}
v_{i+1}=\begin{cases} Fv_i&\text{ if } b_i, b_{i+1}\in B^+\\
Fv_i+\sum_{(d,i+1)\in \mathcal{V}(B)}[a_{d,i+1}]v_d&\text{ if }
b_i\in B^+, b_{i+1}\in B^-\\
V^{-1}v_i&\text{ if } b_i\in B^-, b_{i+1}\in B^+\\
V^{-1}v_i+\sum_{(d,i+1)\in\mathcal{V}(B)}[a_{d,i+1}]v_d&\text{ if }
b_i,b_{i+1}\in B^-.
\end{cases}
\end{equation}
We set 
$$v_i=\sum_{j=0}^{h-1}c_{i,j}e_{b_i+j}$$ with $c_{i,j}\in W(R)$ and write $c_{i,j}=(c_{i,j,l})_{l\in \mathbb{N}}$. Let 
\begin{equation}\label{glphi}\varphi(j,l)=j+lh.
\end{equation}
As $\varphi$ is a bijection between $\{0,\dotsc,h-1\}\times\mathbb{N}$ and $\mathbb{N}$, we may write $\tilde{c}_{i,\varphi(j,l)}$ instead of $c_{i,j,l}$.
Let $\tilde{c}_{i,0}=1$ for all $i$. Then in every point of $S$, the first index of $v_i$ is $b_i$, and its first coefficient is $1$. We define the $\tilde{c}_{i,\varphi}$ by induction on $\varphi(j,l)$, and for fixed $\varphi$ by induction on $i$. For $\varphi>0$ let $\tilde{c}_{0,\varphi}=0$ to satisfy (\ref{glstart}). If $b_{i+1}\in B^+$, (\ref{gldefdisp}) implies that $\tilde{c}_{i+1,\varphi}=\tilde{c}_{i,\varphi}^{\sigma}$. If $b_{i+1}\in B^-$, then $\tilde{c}_{i+1,\varphi}$ is equal to $\tilde{c}_{i,\varphi}^{\sigma}$ plus a polynomial in the $\tilde{c}_{i',\varphi'}$ with $\varphi'<\varphi$ and coefficients in $R$. Hence we can inductively define $\tilde{c}_{i,\varphi(j,l)}=c_{i,j,l}\in R$, and obtain $c_{i,j}=(c_{i,j,l})\in W(R)$.

With these $v_i$ let
\begin{eqnarray*}
L&=&\langle v_i \mid b_i\in B^-\rangle_{W(R)}\\
T&=&\langle v_i \mid b_i\in B^+\rangle_{W(R)}\\
P&=&L\oplus T\\
Q&=&L\oplus I_R T.
\end{eqnarray*} The first indices $b_i$ of the $v_i$ are pairwise non-congruent modulo $h$, hence the $v_i$ are linearly independent over $W(R)[1/p]$ and $P\otimes \mathbb{Q}=N_R$. To show that this defines a display
over $S$, we have to verify that $F(P)\subseteq P$ and that
$V^{-1}(Q)$ generates $P$. The only assertions that do not immediately follow from (\ref{gldefdisp}) are $Fv_{h-1}\in P$ and $v_0\in \langle V^{-1}(Q)\rangle_{W(R)}$. As $l\in A$ for all $l\geq
b_0$, all those $e_l$ are in $P$. As $b_{h-1}+m=b_0$, we have
$Fv_{h-1}=\sum_{j=0}^{h-1} c_{h-1,j}^{\sigma}e_{b_0+j}$ with $w_0(c_{h-1,0})=1$. Therefore, $Fv_{h-1}\in P$, and $Fv_{h-1}=\sum_{i\in B}\delta_iv_i$ with $\delta_i\in W(R)$ and $w_0(\delta_0)=1$. All $v_i$ with $i>0$ are in $\langle V^{-1}(Q)\rangle_{W(R)}$ and $Fv_{h-1}=V^{-1}(pv_{h-1})$. Thus also $v_0\in\langle V^{-1}(Q)\rangle_{W(R)}$.

Let $B$ be the cycle corresponding to a normalised semimodule
$A$. Then this display induces a quasi-isogeny and thus a morphism
\begin{equation}
f_A:\mathbb{A}^{\mathcal{V}(B)}\rightarrow \mathcal{M}_A.
\end{equation}
\begin{lemma}\label{lemimp}
Let $R$ be an excellent local ring and a $k$-algebra. Let $x\in\mathcal{M}_A(R)$ and let $\mathcal{P}=(P,Q,F,V^{-1})$ be the corresponding display over $R$. Assume that for every $b_i\in B^+$ there is a $w_i^0\in P$ and for every $b_i\in B^-$ a $w_i^0\in Q$  with $w_i^0=\sum_{j=0}^{h-1} c_{i,j}^0e_{b_i+j}$ with $c_{i,j}^0\in W(R)$ and the following additional properties. $w_0(c_{i,0}^0)=1$ for all $i$ and the $w_i^0$ generate $P$. Then there is a unique $\tilde{x}\in\mathbb{A}^{\mathcal{V}(B)}(R)$ with $f_A(\tilde{x})=x$.
\end{lemma}
\begin{proof}
We want to show that there exist unique $a_{d,i}\in R$ such that for the corresponding display $\mathcal{P}'=(P',Q',F,V^{-1})$ we have that $P'\subseteq P$ and $Q'\subseteq Q$. As $x$ and $f_A((a_{d,i}))$ are in the same connected component of $\mathcal{M}$, the displays then have to be equal. 

We show by induction on $r$ that for each $b_i\in B^-$
there is a $w_i^{r}\in Q$ and for each $b_i\in B^+$ a $w_i^{r}\in
P$ of the form $w_i^{r}=\sum_{j=0}^{h-1}c_{i,j}^re_{b_i+j}$ with $c_{i,j}^r=(c_{i,j,l}^r)\in W(R)$ and the following property. The coefficients
$c_{i,j,l}^r$ for $\varphi(j,l)\leq r$ are equal to those of the
basis $v_i$ of the display of a point of
$\mathbb{A}^{\mathcal{V}(B)}$ which only depends on $\mathcal{P}$ and not on the chosen basis $w_i^0$. The coordinate $a_{d,i}$ of the point
of $\mathbb{A}^{\mathcal{V}(B)}$ will be determined in the step where
$r=b_d-b_i$. Especially, the point is fixed after finitely many
steps.

For $r=0$ the claim follows from the assumptions of the lemma. As
$b_0+\mathbb{N}\subseteq n+ A$, we can choose $w_0^r=e_{b_0}=v_0\in Q$ for all $r$. Now suppose that the $w_i^{r}$ are defined for some fixed
$r$. We use a second induction on $i$ to define
\begin{equation*}
\tilde{w}_{i+1}^{r+1}=\begin{cases} Fw_i^{r+1}&\text{ if } b_i\in B^+\\
V^{-1}w_i^{r+1}&\text{ if } b_i\in B^-.
\end{cases}
\end{equation*} From $w_i^{r+1}\in Q$ for $b_i\in B^-$ we obtain that $\tilde{w}_{i+1}^{r+1}\in P$ also for those $i$. If $b_{i+1}\in B^+$ let
\begin{equation*}
w_{i+1}^{r+1}=\tilde{w}_{i+1}^{r+1}.
\end{equation*} 
If $b_{i+1}\in B^-$, we now show how to modify $\tilde{w}_{i+1}^{r+1}\in P$ to obtain an element of $Q$. A basis of the free $R$-module $P/Q\cong T/I_RT$ is given by the $w_d^r$ with $b_d\in B^+$. There is an element of $Q$ with first index $b_{i+1}$. Thus there are unique $\alpha_{d,i+1}^{r+1}\in R$ with
\begin{equation}
w_{i+1}^{r+1}=\tilde{w}_{i+1}^{r+1}+\sum_{b_d>b_{i+1}, b_d\in
B^+}[\alpha_{d,i+1}^{r+1}]w_d^r\in Q.
\end{equation}
By the induction hypothesis, the coefficients $c_{i,j,l}^{r+1}$ of $w_i^{r+1}$ with $\varphi(j,l)\leq r+1$ and the coefficients of all $w_d^r$ with $\varphi(j,l)\leq r$ are uniquely defined by $\mathcal{P}$ and independent of the chosen $w_i^0$. This implies $c_{i,j,l}^{r+1}=c_{i,j,l}^r$ for all $j,l$ with $\varphi(j,l)\leq r$. Especially, $\alpha_{d,i+1}^{r+1}=a_{d,i+1}$ for all $b_d\leq b_{i+1}+r$. If $(d,i+1)\in \mathcal{V}(B)$ with
$b_d-b_{i+1}=r+1$, then let $a_{d,i+1}=\alpha_{d,i+1}^{r+1}$. Then $a_{d,i+1}$ also only depends on $\mathcal{P}$. This defines
$w_i^{r+1}$ for all $i$ and unique $a_{d,i}$ satisfying the condition above for $r+1$. 

Each coefficient $c_{i,j,l}^r$ remains fixed after $\varphi(j,l)$ steps. Hence the sequences $w_i^r$ converge in $P$, and their limits $w_i$ are as desired. 
\end{proof}
\begin{theorem}\label{thm43}
Let $A$ be a normalised semimodule. Then
$f_A:\mathbb{A}^{\mathcal{V}(B)}\rightarrow \mathcal{M}_A$ is an
isomorphism.
\end{theorem}
\begin{proof}
If $K\supseteq k$ is a perfect field, each lattice corresponding to an element of $\mathcal{M}_A(K)$ has a basis satisfying the assumptions of Lemma \ref{lemimp}. Hence $f_A(K):\mathbb{A}^{\mathcal{V}(B)}(K)\rightarrow \mathcal{M}_A(K)$ is a bijection. 

We want to show that $f_A$ is proper by verifying the valuation
criterion. Let $x\in \mathcal{M}_A(k[[t]])$ and let $x_{\eta}$ and
$x_0$ be its generic and special point. Let $\tilde{x}_{\eta}\in
\mathbb{A}^{\mathcal{V}(B)}(k((t)))$ be a point mapping to
$x_{\eta}$. Let $\mathcal{P}=(P,Q,F,V^{-1})$ be the display of $x$.
The $W(k[[t]])$-module $P$
is a submodule of $P_{\eta}=P\otimes_{W(k[[t]])}W(k((t)))$, the first component of the
display $\mathcal{P}_{\eta}$ of $x_{\eta}$. As $\tilde{x}_{\eta}$ maps to $x_{\eta}$,
we can describe $\mathcal{P}_{\eta}$ as generated by elements $v_i$ as
above. By Lemma \ref{lemimp}, $x_0\in \mathcal{M}_A(k)$ is also in the image of $f_A$. Hence we can choose generators $v_i'$ of $\mathcal{P}$ which for $b_i\in B^-$ are in $Q$, and which modulo $(t)$ reduce
to the standard generators of the display of the inverse image of
$x_0$ under $f_A$. Let $y$ be the minimal element of $B$. As $v_i'\in N_{k[[t]]}$, there is an $a\in \mathbb{N}$ such that $v_i'\in \langle e_y,\dotsc,e_{y+h-1} \rangle_{W(k[[s]])}$ where $s^{p^{a}}=t$. In the following we consider $x$ as a $k[[s]]$-valued point of $\mathcal{M}_{A}$. The reduction of $v_i'$ modulo $(s)$ has first index $b_i$. As $x_{\eta}\in \mathcal{M}_A$, the index of the first nonzero coefficient of each $v_i'$ is in $A$. Thus we can modify each $v_i'$ by a
linear combination of the $v_j'$ with $b_j<b_i$ and coefficients
in $W(k[[s]])$ which reduce to $0$ modulo $(s)$ such that the new first nonzero coefficient is that
of $e_{b_i}$. Besides, we only have to modify the elements $v_i\in Q$ by other elements of $Q$. Therefore we may in addition assume that the first
nonzero coefficient of $v_i'$ has index $b_i$, and is $1$. After passing to a larger $a$ we may assume that $v_i'\in \langle e_{b_i},\dotsc,e_{b_i+h-1}\rangle_{W(k[[s]])}$. By Lemma \ref{lemimp} for $w_{i}^0=v_i'$ we obtain a unique point $\tilde{x}'=(b_{d,i})\in\mathbb{A}^{\mathcal{V}(B)}(k[[s]])$ mapping to $x$. But as its generic point $\tilde{x}'_{\eta}$ maps to $x_{\eta}$, the uniqueness in Lemma \ref{lemimp} implies that $\tilde{x}'_{\eta}=\tilde{x}_{\eta}\in \mathbb{A}^{\mathcal{V}(B)}(k((t)))$. Hence $b_{d,i}\in k[[s]]\cap k((t))=k[[t]]$, and $x$ is in the image of $f_A(k[[t]])$.

Lemma \ref{lemimp} further implies that the tangent morphism of $f_A$ is injective at every closed point. The theorem now follows since a morphism of reduced schemes of finite type over an algebraically closed field, which is finite, universally bijective, and whose tangent morphism is injective at every closed point is an isomorphism (compare \cite{Inken}, Lemma 5.13).
\end{proof}

Using the paving of $\mathcal{M}_{\red}^0$ by affine spaces and
that $\mathcal{M}_{\red}^0$ is projective, one obtains the
following result about its cohomology.
\begin{theorem}\label{thmcohom}
Let $l\neq p$ be prime. Then $H^{2i+1}(\mathcal{M}_{\red}^0,\mathbb{Q}_l)=0$ and  $H^{2i}(\mathcal{M}_{\red}^0,\mathbb{Q}_l)$ is a successive extension of $d(i)$ copies of $\mathbb{Q}_l(-i)$ for all $i$. Here $d(i)$ is the number of normalised cycles $B$ with $\mid \mathcal{V}(B)\mid =i$.
\end{theorem}

\begin{proposition}\label{kordimkoh}
\begin{mylist}
\item $d(0)=d((m-1)(n-1)/2)=1$ for all $m$ and $n$. If
$m,n>1$, also $d(1)=1$.
\item Let $\min\{m,n\}=2$. Then $d(i)=1$ for $0\leq i\leq \dim \mathcal{M}_{\red}^0$.
\item Let $\min\{m,n\}>2$. Then $d((m-1)(n-1)/2-1)>1$.
\end{mylist}
\end{proposition}
\begin{proof}
The equation $d((m-1)(n-1)/2)=1$ is shown in
\cite{deJongOort}, 6. They also show that for a semimodule $A$,
the dimension $ \mid \mathcal{V}(B(A)) \mid $ is bounded below by the number
of positive integers $s$ such that there exists an $a\in A$ with
$a+s\notin A$ (see \cite{deJongOort}, 6.12). Let $A$ be a
normalised semimodule with $ \mid \mathcal{V}(B(A)) \mid =0$. Then $a\in A$
implies $a'\in A$ for all $a'>a$. Thus $A=\mathbb{N}$. One easily
sees that this semimodule indeed leads to a zero-dimensional
subscheme. Let now $A$ be normalised with $ \mid \mathcal{V}(B(A)) \mid =1$.
Then $a\in A$ implies that there is at most one element of
$\mathbb{Z}\setminus A$ that is larger than $a$. Analogously, for
$a\notin A$ there is at most one element of $A$ smaller than $a$.
This leaves only $A=\{-1,1,2,3,\dotsc\}$ as a candidate for a
contribution to $d(1)$. It is a semimodule if and only if
$m,n>1$. Again one can see (using the combinatorics
explained in \cite{deJongOort}, 6) that $ \mid \mathcal{V}(B(A)) \mid =1$.

To show (ii), we may assume that $m=2$ and $n=2l+1$ for some $l$.
Each normalised semimodule is of the form
$A=A_i=(2\mathbb{N}-i)\cup(\mathbb{N}+i+1)$ for some $i\in
\{0,\dotsc,l\}$.  The cycle $B_i=B(A_i)$ is
$$(2l+2+i,i+1,i+3,\dotsc,2l-i-1,2l-i+1,-i,-i+2,\dotsc,2l+i)$$ with
$B^-=\{2l+2+i,2l-i+1\}$. The element $2l+2+i$ is the largest
element of $B_i$, so it does not contribute to $\mathcal{V}(B_i)$.
The other element of $B^-$ is smaller than the $i$ elements
$2l-i+2,\dotsc,2l+i$ of $B^+$. Hence $ \mid \mathcal{V}(B_i) \mid =i$.

For (iii) we have to construct two normalised semimodules
leading to subschemes of codimension $1$. Assume that $m<n$, the
other case is completely analogous. Let
\begin{equation}
A_1=\{am+bn \mid a, b\geq 0\}\cup \{mn-m-n\}.
\end{equation} There are $(m-1)(n-1)/2$ natural numbers
which cannot be written as $am+bn$ with $a,b\geq 0$, and $mn-m-n$
is the largest. Thus the lower bound on $ \mid \mathcal{V}(A_1)
\mid $ used in the proof of (i) shows that the codimension of the
subscheme corresponding to the normalisation of $A_1$ is at most
$1$. But there is only one semimodule leading to a subscheme of
codimension $0$, and this is obtained by normalising $A_0=\{am+bn \mid
a,b\geq 0\}$ (see \cite{deJongOort}, 6). Thus the
normalisation of the semimodule $A_1$ leads to a subscheme of
codimension $1$. We now define a second subscheme of codimension 1. The cycle corresponding to $A_0$ is given by
$B^+=\{0,m,\dotsc,(n-1)m\}$ and $B^-=\{n,2n,\dotsc,mn\}$. Let
$B_2^+$ be the index set obtained from $B^+$ by replacing $(n-1)m$
by $(n-1)m-2n$, and let $B_2^-$ be obtained from $B^-$ by
replacing $mn$ by $mn-m$ and $mn-n$ by $mn-m-n$. One can easily check that this defines a cycle. A pair of
elements $(i,j)\in B^+\times B^-$ with $i>j$ is then replaced by a
pair in $B^+_2\times B^-_2$. The pair $((n-1)m,mn-n)$ which is
replaced by $(mn-m-2n,mn-m-n)$ is the only pair with larger first entry which is replaced
by a pair such that the first entry is smaller than the second.
After normalising the cycle we get again a subscheme of codimension $1$. The smallest element of $A_1$ and $A_2$ is 0. As $mn-m-2n\in A_2\setminus A_1$, the normalisations of the two semimodules are different. 
\end{proof}

\subsection{Application to smoothness}

In this section we show the following
\begin{theorem}\label{thmsmooth} Let $\mathbb{X}$ be an arbitrary
$p$-divisible group over an algebraically closed field of
characteristic $p$. Then $\mathcal{M}_{\red}^0$ is smooth if and
only if one of the following holds. Either $\dim\mathcal{M}_{\red}^0=0$
or the isocrystal $N$ of $\mathbb{X}_{\bi}$ is simple of slope
$2/5$ or $3/5$.
\end{theorem}
\begin{remark}
The condition $\dim\mathcal{M}_{\red}^0=0$ is equivalent to the
condition that $\mathbb{X}$ is ordinary or that the isocrystal of $\mathbb{X}_{\bi}$ is simple of
slope $m/(m+n)$ with $\min\{m,n\}=1$.
\end{remark}
Once we have shown the theorem for bi-infinitesimal $\mathbb{X}$, we can treat the general case as in Section
\ref{secmultet}. We may thus assume that $\mathbb{X}$ is
bi-infinitesimal. The results of Sections \ref{secconncomp} and
\ref{secirrcomp} imply that the connected components of
$\mathcal{M}_{\red}$ are irreducible if and only if $N$ is simple.
From now on we assume this. Let $m/(m+n)$ be its slope with
$(m,n)=1$. We also assume that $\id_{\mathbb{X}}$ corresponds to a
lattice of volume $0$. We consider the
following cases.\\

{\bf Case 1:} $\min\{m,n\}=1.$\\

\noindent If $m$ or $n$ is 1, the dimension of $\mathcal{M}_{\red}^0$ is $0$ and the scheme is smooth.\\

{\bf Case 2: $\{m,n\}=\{2,3\}.$}

\begin{theorem}
Let $\mathbb{X}$ be bi-infinitesimal and let its rational Dieudonn\'{e}
module $N$ be simple of slope $2/5$ or $3/5$. Then
$\mathcal{M}_{\red}^0\cong \mathbb{P}^1$.
\end{theorem}

\begin{proof}We assume that $m=2$ and $n=3$. The case $n=2$ and $m=3$
is similar and thus omitted.

Let $\mathbb{P}^1=U_0\cup U_1$ be the standard open covering. We
denote the points of $\mathbb{P}^1$ by $[a_{-1}:a_0]$. Over $U_0\cong \Spec(k[a_0])$
let
\begin{eqnarray*}
L_0&=&\langle e_{2}+[a_0]e_3,e_5\rangle_{W(k[a_0])}\\
T_0&=&\langle e_{-1}+[a_0]^{\sigma}e_0,e_1,e_3\rangle_{W(k[a_0])}
\end{eqnarray*}
and
\begin{eqnarray*}
P_0&=&L_0\oplus T_0=\langle e_{-1}+[a_0]^{\sigma}e_0,e_1,e_2,e_3,e_5 \rangle_{W(k[a_0])}\\
Q_0&=&L_0\oplus I_{k[a_0]}T_0.
\end{eqnarray*}
On the other hand let $P$ and $Q$ be the display from the
definition of $f_A$ for the semimodule $A=\{-1,1,2\dotsc\}$. Then
$\{e_1, e_2,\dotsc\}\subset P$ and $\{e_4,e_5,\dotsc\}\subset Q$.
Using this and the first steps of the recursion for the
generators of $P$, one can see that $P=P_0$ and $Q=Q_0$. Thus $P_0$ and $Q_0$
define a display. As $A$ is the minimal semimodule, the
corresponding morphism $\mathbb{A}^1\rightarrow \mathcal{M}_{\red}^0$ is
an open immersion. Over $U_1\cong \Spec (k[a_{-1}])$ let
\begin{eqnarray*}
L_1&=&\langle [a_{-1}]e_2+e_3,e_4 \rangle_{W(k[a_{-1}])}\\
T_1&=&\langle [a_{-1}]^{\sigma} e_{-1}+e_0,e_1,
e_2\rangle_{W(k[a_{-1}])}
\end{eqnarray*}
and choose
\begin{eqnarray*}
P_1&=&L_1\oplus T_1=\langle [a_{-1}]^{\sigma}e_{-1}+e_0,e_1,e_2,e_3,e_4 \rangle_{W(k[a_{-1}])}\\
Q_1&=&L_1\oplus I_{k[a_{-1}]}T_1.
\end{eqnarray*}
One easily checks that this defines a display. It
is obvious that the corresponding morphism
$\varphi_1:\mathbb{A}^1\rightarrow \mathcal{M}_{\red}^0$ is injective on
$R$-valued points. As for $f_A$ in Theorem \ref{thm43} one can show that $\varphi_1$ is an immersion. The complement of its image consists of the image
of the origin in $U_0$. We can glue the
morphisms corresponding to the displays over $U_0$ and $U_1$ to obtain an isomorphism $\mathbb{P}^1\rightarrow \mathcal{M}_{\red}^0$.
\end{proof}

{\bf Case 3: $\min\{m,n\}=2$ and $\max\{m,n\}>3.$}\\

\noindent We consider the case $m=2$ and $n=2l+1$ with $l>1$. The
case $n=2$ and $m=2l+1$ is similar and thus omitted. We have
$\dim\mathcal{M}_{\red}^0=l$. Let $M\subset N$ be the
Dieudonn\'{e} lattice generated by $e_{-l+2}$ and $e_{l-1}$. Then
this lattice corresponds to the $k$-valued point $x=f_{A_{l-2}}(0)\in\mathcal{M}_{A_{l-2}}(k)$ where $A_{l-2}$ is as in the proof of
Proposition \ref{kordimkoh}.
\begin{proposition}
The dimension of the tangent space of $\mathcal{M}_{\red}^0$ in $x$ is at
least $l+1$.
\end{proposition}
\begin{proof} For
$(a_0,\dotsc,a_{l-1},b_0,b_1)\in k^{l+2}$ consider the following
submodules of $N_{k[\varepsilon]}$ where $k[\varepsilon]\cong
k[t]/(t^2)$. Let
\begin{equation}
v_1=e_{l+3}+[\varepsilon]([a_0]e_{l+1}+\sum_{i=1}^{l-1}[a_i]e_{l+2i})
\end{equation}
and 
\begin{equation}\label{gldefv2e}
v_2=e_{3l}+[\varepsilon b_0]e_{l+1}+[\varepsilon b_1]e_{l+2}.
\end{equation} 
Let
\begin{align*}
L&=\langle v_1,v_2 \rangle_{W(k[\varepsilon])}\\
T&=\langle
e_{-l+2},e_{-l+4},\dotsc,e_{3l-2},e_{l-1},e_{l+1}\rangle_{W(k[\varepsilon])},
\end{align*}
then
\begin{align}
P&=L\oplus T=M\otimes_{W(k)}W(k[\varepsilon])\\
\label{glqtang}Q&=L\oplus I_{k[\varepsilon]}T.
\end{align}
As $\sigma(\varepsilon)=0$, this defines a display. For $i_0\in \{l+1,l+2,l+4,\dotsc,3l-2\}$ there is no element of $Q$ of the form $\sum_{i\geq i_0}[\delta_i ]e_i$ with $\delta_i\in k[\varepsilon]$ and $\delta_{i_0}\neq 0$. This implies that the display leads to an $l+2$-dimensional
subspace of the tangent space of $\mathcal{M}$ at $x$. We now have
to construct an $l+1$-dimensional subspace that lies in the tangent space of $\mathcal{M}_{\red}$. For $a_0,a_1\neq 0$ let 
\begin{equation}\label{gldefv1}
v_1=[ta_0]e_{l+1}+[ta_1]e_{l+2}+e_{l+3}+\sum_{i=2}^{l-1}[ta_i]e_{l+2i}\in N_{k((t))}.
\end{equation} 
Let further 
\begin{align*}
L_1&=\langle e_{3l+2}, v_1\rangle_{W(k((t)))}\\
T_1&=\langle V^{-1}v_1, FV^{-1}v_1,\dotsc,F^{2l}V^{-1}v_1\rangle_{W(k((t)))}\\
P_1&=L_1\oplus T_1\\
Q_1&=L_1\oplus I_{k((t))}T_1.
\end{align*}
As $a_0\neq 0$, there is an element of $P_1$ with first index $i$ for all $i\geq l$ and of $Q_1$ for all $i\geq 3l+1$. Using this one can easily see that $F(P_1)\subseteq P_1$ and that $V^{-1}(Q_1)$ generates $P_1$. Thus $P_1$ and $Q_1$ define a display over $k((t))$, and a $k((t))$-valued point of $\mathcal{M}^0_{\red}$. As $\mathcal{M}_{\red}^0$ is projective this point is induced by a $k[[t]]$-valued point. Its display is $(P_1\cap  N_{k[[t]]}, Q_1\cap  N_{k[[t]]}, F, V^{-1})$. We want to show that the special point corresponds to $x$. The element $V^{-1}v_1$ of $P_1\cap N_{k[[t]]}$ reduces to $V^{-1}(e_{l+3})=e_{-l+2}$ modulo $(t)$. To show that $M$ is contained in the reduction of $P_1\cap N_{k[[t]]}$ modulo $(t)$, it remains to see that $e_{l-1}$ is contained in this reduction. For all $i\geq 3l+1$ the vector $e_i$ is in $Q_1$. We consider the following element of $Q_1$ modulo the lattice generated by these elements $e_i$.
\begin{align*}
(F^{l-1}-[(ta_0)^{\sigma^{l-1}}]F^{l-2})v_1&\equiv
[t^{\sigma^{l-1}}]\left( -[a_0^{\sigma^{l-1}}(ta_0)^{\sigma^{l-2}}]e_{3l-3}
-[a_0^{\sigma^{l-1}}(ta_1)^{\sigma^{l-2}}]e_{3l-2}\right.\\
&\hspace{3cm}\left.+([a_1^{\sigma^{l-1}}]-[a_0^{\sigma^{l-1}}(ta_2)^{\sigma^{l-2}}])e_{3l} \right)\\
&= [t^{\sigma^{l-1}}]v
\end{align*}
for some $v\in Q_1\cap N_{k[[t]]}$. The reduction of $v$ modulo $(t)$ is $[a_1^{\sigma^{l-1}}]e_{3l}$. Thus $e_{l-1}=V^{-1}(e_{3l})$ is contained in the lattice at the special point. Hence the special point of this $k[[t]]$-valued point is $x$. If $l>2$, the reduction of $v$ modulo $(t^2)$ is $[a_1^{\sigma^{l-1}}]e_{3l}$. Hence $e_{3l}$ is in the projection of $Q_1\cap N_{k[[t]]}$ to $N_{k[\varepsilon]}$. If $l=2$, the reduction of $([a_1^{\sigma^{l-1}}]-[a_0^{\sigma^{l-1}}(ta_2)^{\sigma^{l-2}}])^{-1}v$ modulo $(t^2)$ is equal to $v_2$ as in (\ref{gldefv2e}) with $b_0=-a_0^{\sigma}a_0/a_1^{\sigma}$ and $b_1=-a_0^{\sigma}a_1/a_1^{\sigma}$. Comparing the image of $Q_1\cap N_{k[[t]]}$ under the projection to $N_{k[\varepsilon]}$ to the definition of $Q$ in (\ref{glqtang}) we see that the tangent vector of this $k[[t]]$-valued point at $x$ corresponds to the tangent vector $(a_0,\dotsc,a_{l-1},0,0)\in k^{l+2}$ if $l>2$ and to $$(a_0,a_{1}, \frac{-a_0^{\sigma}a_0}{a_1^{\sigma}},\frac{-a_0^{\sigma}a_1}{a_1^{\sigma}})$$ if $l=2$.

For $b_0\neq 0$ let $$v_2=[tb_0]e_{l+1}+e_{3l}$$ and  
\begin{align*}
L_2&=\langle e_{3l+2}, v_2\rangle_{W(k((t)))}\\
T_2&=\langle V^{-1}v_2, FV^{-1}v_2,\dotsc,F^{2l}V^{-1}v_2\rangle_{W(k((t)))}\\
P_2&=L_2\oplus T_2\\
Q_2&=L_2\oplus I_{k((t))}T_2.
\end{align*}
The same reasoning as above shows that this defines a display over $k((t))$. To show that it leads to a $k[[t]]$-valued point of $\mathcal{M}_{\red}$ with special point $x$, we have to check that $e_{l-1}$ and $e_{-l+2}$ are in the lattice at the special point. The reduction of $V^{-1}v_2=[tb_0]^{\sigma}e_{-l}+e_{l-1}$ modulo $(t)$ is $e_{l-1}$. Besides, $Fv_2-e_{3l+2}=[t^{\sigma}b_0^{\sigma}]e_{l+3}\in Q_2$, hence $e_{l+3}\in Q_2\cap N_{k[[t]]}$. As $e_{-l+2}=V^{-1}e_{l+3}$, the lattice $M$ is contained in the reduction of $P_2$ modulo $(t)$. The fact that $e_{l+3}\in Q_2\cap N_{k[[t]]}$ also shows that the tangent vector of this $k[[t]]$-valued point in $x$ corresponds to $(0,\dotsc,0,b_0,0)\in k^{l+2}$. Thus we constructed elements of the tangent space of $\mathcal{M}_{\red}$ in $x$ which generate an $l+1$-dimensional subspace.
\end{proof}

{\bf Case 4: $\min\{m,n\}>2$}\\

\noindent In this case Theorem \ref{thmcohom} and Proposition \ref{kordimkoh} show that $$\dim H^2(\mathcal{M}_{\red}^0, \mathbb{Q}_l)\neq\dim H^{2\dim{\mathcal{M}_{\red}^0}-2}(\mathcal{M}_{\red}^0, \mathbb{Q}_l).$$ Hence $\mathcal{M}_{\red}^0$ does not satisfy Poincar\'{e} duality and it cannot be smooth.
\qed

\end{document}